\documentclass[11pt]{amsart}
\usepackage[T1]{fontenc}
\usepackage{lmodern,amsmath,amssymb,mathtools,mathrsfs,booktabs}
\usepackage{xcolor}
\usepackage{tikz}
\usetikzlibrary{arrows.meta}
\usepackage[margin=1.1in]{geometry}
\usepackage[colorlinks=true,linkcolor=blue,citecolor=blue,urlcolor=blue]{hyperref}
\newcommand{\rev}[1]{\leavevmode{#1}}
\newtheorem{theorem}{Theorem}[section]
\newtheorem*{introtheorem}{Theorem}
\newtheorem{proposition}[theorem]{Proposition}
\newtheorem{lemma}[theorem]{Lemma}
\newtheorem{corollary}[theorem]{Corollary}
\theoremstyle{definition}

\newcommand{\C}{\mathbb C}
\newcommand{\PP}{\mathbb P}
\newcommand{\M}{\overline{\mathcal M}}
\newcommand{\OO}{\mathcal O}
\newcommand{\Res}{\operatorname{Res}}
\newcommand{\diag}{\operatorname{diag}}
\newcommand{\End}{\operatorname{End}}
\newcommand{\id}{\mathrm{id}}
\newcommand{\pt}{\mathrm{pt}}
\newcommand{\bq}{\mathbf q}
\newcommand{\xh}{\hat x}
\newcommand{\yt}{\widetilde y}
\newcommand{\eh}{\widehat e}
\newcommand{\thh}{\widehat\theta}
\newcommand{\GW}{\mathrm{GW}}
\newcommand{\TR}{\mathrm{TR}}
\DeclareMathOperator{\st}{st}
\DeclareMathOperator{\ev}{ev}
\title[All-genus mirror symmetry for type D]{All-genus mirror symmetry for the type-$D$ minimal resolution Calabi-Yau threefold}
\author{Bohan Fang}
\address{Beijing International Center for Mathematical Research, Peking University, 5 Yiheyuan Road, Beijing 100871, China}
\email{bohanfang@gmail.com}
\author{Zhuoming Lan}
\address{Beijing International Center for Mathematical Research, Peking University, 5 Yiheyuan Road, Beijing 100871, China}
\email{lanzm26@pku.edu.cn}
\date{}
\begin{document}
\begin{abstract}
Let $Y$ be the product of a minimal resolution of a binary dihedral surface singularity with an equivariant line, with volume-preserving scalar action. We identify its stable equivariant descendant Gromov--Witten theory and its free energies in genus at least two with the sign-sector topological recursion of the degenerate type-$D$ relativistic Toda spectral curve. The $R-$matrix calibration is determined at a partial large-radius limit, where the trivalent exceptional curve retains a nonzero K\"ahler parameter, while the remaining parameters vanish. The limiting theory on the mirror side is a sum of twisted point theories and the local $(-2,0)$ curve theory. Moreover, we identify the equivariant K-group supported on the surface zero section with an integral relative homology group of the mirror. Explicit cycles for the dual tautological bundles realize their Gamma central charges and all stable Gamma-framed descendants.
\end{abstract}
\maketitle

\section{Introduction}
Let $G\subset SU(2)$ be the binary dihedral group corresponding to $D_l$, where $l\geq4$, and let $Z\to\C^2/G$ be its minimal resolution. We consider
\[
 Y=Z\times\C_{-2\nu}.
\]
The torus $T=\C^*$ acts on $Z$ through scalar multiplication of weight $\nu$ on $\C^2$. Together with the weight $-2\nu$ on the extra line, this gives a Calabi--Yau action. We localize equivariant invariants over $\C(\nu)$.

The remodeling conjecture of Bouchard--Klemm--Mari\~no--Pasquetti \cite{BKMP}, building on Eynard--Orantin recursion \cite{EO} and Mari\~no's work on open topological strings \cite{Marino}, relates all-genus toric Calabi--Yau Gromov--Witten theory to recursion on mirror curves. It was proved for toric threefolds by Eynard--Orantin \cite{EOBKMP} and for toric three-orbifolds by Fang--Liu--Zong \cite{FLZ}. The graph comparison connects Givental's higher-genus formalism and Teleman's reconstruction \cite{GiventalHigher,GiventalQuant,Teleman} with the intersection-theoretic description of recursion due to Eynard and Dunin-Barkowski--Orantin--Shadrin--Spitz \cite{Eynard,DOSS}.

For ADE resolutions, Bryan--Gholampour \cite{BG} express the equivariant quantum product in terms of the root system. The type-$D_l$ affine relativistic Toda spectral curve was explicitly constructed by Borot--Brini \cite[Section~6.3]{BB}. The mirror curve used here arises from its degeneration at vanishing base Novikov parameter. Brini discusses the geometric interpretation of this degeneration as a mirror of an ADE surface times a line, with an explicit analysis in type $E_8$ \cite[Section~4.3.2]{BriniE8}. The genus-zero mirror theorem of Brini--Ma--Strachan \cite[Theorem~3.2]{BMS} identifies the resulting Landau--Ginzburg Frobenius structure with the equivariant quantum cohomology of $Z$.

\rev{A recent work of Brini--Ma establishes the integral and higher-genus quantum McKay correspondence in type $D$, including the Calabi--Yau threefold extension \cite[Theorems~5.12--5.13 and Corollaries~6.12 and~7.1]{BM}.} Fang--Lan--Ma \cite{FLM} construct the sign-sector recursion and identify it with the equivariant Gromov--Witten theory of $[\C^2/G]\times\C$. \rev{In this paper, we establish the resolution-side remodeling comparison directly from a partial large-radius boundary of the K\"ahler moduli of $Y$.}

The scalar action fixes the trivalent exceptional curve $E$ pointwise. The isolated fixed points of $Y$ are $p_1,\ldots,p_{l-3},p_L,p_R$, defined geometrically in Section~\ref{sec:amodel}. The classical cohomology of $E$ has a nilpotent divisor class, so the equivariant quantum algebra at the full large-radius limit is not semisimple. Retaining the parameter $Q$ of $E$ gives a semisimple boundary theory. When all other curve parameters vanish, the localized cohomological field theory splits as
\begin{equation}\label{eq:intro-splitting}
 \Omega^Y\big|_{\bq=0}
 =\Omega^{\PP^1,\,\OO(-2)_{2\nu}\oplus\OO_{-2\nu}}
  \oplus\bigoplus_{a\in\{1,\ldots,l-3,L,R\}}\Omega^{\pt,T_{p_a}Y}.
\end{equation}
The Toda curve has the same decomposition in its local critical-point data. Beta integrals compute the point calibrations, and toric remodeling \cite{FLZ} determines the rank-two curve calibration. This fixes the constant diagonal ambiguity left by genus-zero mirror symmetry.

\rev{We developed our partial large-radius boundary method independently. A corresponding boundary is used in Brini--Ma \cite[Section~6.3]{BM} for the $R$-matrix comparison. Here the boundary determines both the $R$-matrix calibration and the Gamma normalization of the correspondence between supported equivariant K-classes and integral mirror cycles.}

The following statement combines Proposition~\ref{prop:r-equality} and Theorems~\ref{thm:descendants} and~\ref{thm:free-energies}. The leaf substitution is defined in Section~\ref{sec:correlators}.
\begin{introtheorem}
Let $R^Y$ be the canonical calibration selected by the localized equivariant Gromov--Witten CohFT of $Y$, and let $R^B$ be the stationary-phase calibration of the reduced Toda curve with kernel $B^-$. Under the genus-zero mirror map, with compatible choices of critical points and normalized idempotents,
\[
 R^Y(z)=R^B(z).
\]
For $n>0$ and $2g-2+n>0$, the recursion differential $\omega_{g,n}$ gives the equivariant descendant generating function of $Y$ under the substitution \eqref{eq:descendant-dictionary}. For $g\geq2$,
\begin{equation}\label{eq:main-free}
 F_g^{\GW}(Y)=F_g^{\TR}
 =\frac{1}{2-2g}\sum_\alpha
   \Res_{p=p_\alpha}\Phi_\alpha(p)\omega_{g,1}(p),
 \qquad d\Phi_\alpha=\yt_\alpha\,d\xh.
\end{equation}
These are coefficientwise identities in the large-radius expansion, on a finite cover admitting normalized idempotents. The B-model expressions give their continuation through adjacent semisimple domains.
\end{introtheorem}

The two-point input throughout is
\[
 B^-(\mu_1,\mu_2)=
 \left(\frac{1}{(\mu_1-\mu_2)^2}
       +\frac{1}{(\mu_1\mu_2-1)^2}\right)d\mu_1d\mu_2.
\]
Its standard diagonal coefficient fixes the normalization of the recursion and the oscillatory integrals.

Hosono proposed central-charge formulas relating coherent sheaves to integral periods in local mirror symmetry \cite{Hosono}. Gamma classes enter the rational structure considered by Katzarkov--Kontsevich--Pantev \cite{KKP} and the K-theoretic integral structure of Iritani \cite{Iritani}. Fang realized equivariant central charges of compact semipositive toric orbifolds as oscillatory integrals over mirror cycles \cite{Fang}. The all-genus descendant correspondence was established for the equivariant projective line by Fang--Liu--Zong \cite{FLZP1} and for toric Calabi--Yau three-orbifolds by Fang--Liu--Yu--Zong \cite{FLYZ}. In these results, integrating each external variable along a cycle representing a K-class realizes the descendant substitution.

For the type-$D$ Toda curve, Ma relates twisted periods to the quantum differential equation and studies their large-radius asymptotics \cite{MaThesis}. \rev{Here we construct explicit integral relative cycles for supported K-classes. We determine the resolution-side period normalization directly at a partial large-radius boundary. This provides a different approach compared with \cite[Section~4.3]{BM}, in which the resolution-side period matrix is computed from full large-radius asymptotics.}

For the resolution, the relevant group is $K_Z^T(Y)$, the Grothendieck group of equivariant coherent sheaves supported on the surface zero section. Our second result combines Theorems~\ref{thm:k-cycles}, \ref{thm:gamma-period}, and~\ref{thm:gamma-descendants}. In the following theorem, $\int^{\rm reg}$ denotes meromorphic continuation in $2\nu/z$ from a domain where the contour integral converges, with logarithmic branches fixed.
\begin{introtheorem}
Let $\mathcal H^+$ be the integral sign-twisted relative homology of the quotient Toda curve, with ends allowed at the poles of $\lambda$. There is an isomorphism
\[
 \operatorname{mir}_Z:K_Z^T(Y)\longrightarrow\mathcal H^+,
\]
semilinear under $\mathsf u\mapsto\mathsf q^{-1}$, whose values on the supported line bundles are given in \eqref{eq:line-bundle-cycles}. At $\mathsf q=e^{2\pi i\nu/z}$, the oscillatory integrals along the chosen lifts satisfy
\[
 \eta\bigl(\kappa_z(F),J_Y^{\rm norm}(\tau,-z)\bigr)
 =-\frac{\pi iz}{\nu}
       \int_{\operatorname{mir}_Z(F)}^{\rm reg}e^{-W/z}\,d\log\mu.
\]
For every stable $(g,n)$ with $n>0$, integrating $\omega_{g,n}$ over these cycles gives the equivariant Gamma-framed descendant invariants, with factor $(2\pi)^n$, as in \eqref{eq:gamma-descendants}.
\end{introtheorem}
Here $\mathsf u$ is the character of weight $\nu$ and $\mathsf q$ is the coefficient variable of the local system. The Gamma framing, divisor coordinates, branches, and regularization are specified in Sections~\ref{sec:cycles} and~\ref{sec:gamma-periods}. The proof normalizes the integral contours on the components of the partial large-radius degeneration. Point components give beta integrals, while the trivalent component retains a rank-two quantum equation. These boundary data determine the global solution.

\subsection*{Acknowledgments}
We thank Andrea Brini and Jingxiang Ma for discussions on type-$D$ mirror symmetry and the crepant resolution correspondence.

Bohan Fang and Zhuoming Lan are partially supported by the National Key R\&D Program of China (2023YFA1009803) and NSFC grant 12125101. Bohan Fang is also partially supported by NSFC grant 11890661.

\subsection*{AI Usage}

AI tools assisted with parts of the mathematical work and writing. The authors checked and revised the generated material.

\section{The resolution theory at a partial large-radius limit}
\label{sec:amodel}
Let $G\subset SU(2)$ be the binary dihedral group of type $D_l$, $l\geq4$, and let $Z\to\C^2/G$ be the minimal resolution. Scalar multiplication on $\C^2$ induces an action of $T=\C^*$ on $Z$. Write $\nu\in H_T^2(\pt)$ for the weight of the standard character, and set
\[
 Y=Z\times\C_{-2\nu},
\]
where $\C_{-2\nu}$ is the one-dimensional representation of weight $-2\nu$. The holomorphic symplectic form on $Z$ has weight $2\nu$, so the product action preserves the holomorphic volume form on $Y$.

The exceptional locus of $Z$ consists of $l$ rational curves of self-intersection $-2$, meeting according to the $D_l$ Dynkin diagram. Let $\epsilon_1,\ldots,\epsilon_l$ be the standard orthonormal basis of $\mathbb R^l$. Label the exceptional curves $E_i$ by the simple roots
\[
 \alpha_i=\epsilon_i-\epsilon_{i+1}\quad(1\leq i<l),\qquad
 \alpha_l=\epsilon_{l-1}+\epsilon_l.
\]
Their intersection numbers are $E_i\cdot E_j=-(\alpha_i,\alpha_j)$. Thus $E=E_{l-2}$ is the trivalent curve, meeting $E_{l-3}$, $E_{l-1}$, and $E_l$. Write $q_i$ for the Novikov variable of $E_i$, put $Q=q_{l-2}$, and let $\bq=(q_i)_{i\ne l-2}$. The \emph{partial large-radius boundary} considered below is $\bq=0$ with $Q$ retained. \rev{The corresponding surface boundary and calibration are also considered in \cite[Section~6.3]{BM}.}

We compute the localized CohFT at $\bq=0$ and specialize its reconstruction matrix to this boundary. The retained parameter $Q$ makes the curve block semisimple, as shown in Section~2.2.

\subsection{Fixed geometry and localization}
The three points $E\cap E_j$ for $j=l-3,l-1,l$ are fixed, so the action on $E\simeq\PP^1$ is trivial. Let $p_1$ be the $T$-fixed point on $E_1$ outside $E_2$, and put $p_k=E_{k-1}\cap E_k$ for $2\le k\le l-3$. Let $p_L$ and $p_R$ be the $T$-fixed points on $E_{l-1}$ and $E_l$, respectively, outside $E$. The surface symplectic form has weight $2\nu$, giving
\begin{equation}\label{eq:fixed-geometry}
 Y^T=E\sqcup\{p_1,\ldots,p_{l-3},p_L,p_R\},\qquad
 V:=N_{E/Y}=\OO_E(-2)_{2\nu}\oplus\OO_{E,-2\nu}.
\end{equation}
Here the subscripts record the torus weights on the fibers over the fixed curve $E$. At distance $m$ along an arm, the isolated fixed point has tangent weights
\begin{equation}\label{eq:point-weights}
 (-2m\nu,\,2(m+1)\nu,\,-2\nu),
 \qquad m\geq1.
\end{equation}

The long arm contributes one point for each distance $m=1,\ldots,l-3$, and each short arm contributes a point with $m=1$.

Set $H=H_T^*(Y)\otimes\C(\nu)$. Its dimension is $l+1$, and its pairing is
\[
 \eta(a,b)=\int_{Y^T}\frac{i^*(ab)}{e_T(N_{Y^T/Y})}.
\]
For $2g-2+n>0$, put $\mathcal M_\beta=\M_{g,n}(Y,\beta)$, with universal curve $\pi$ and universal map $f$. The complex $(R\pi_*f^*TY)^\vee$ gives the equivariant perfect obstruction theory relative to the smooth stack of prestable curves. The associated absolute theory defines a virtual class in $A_*^T(\mathcal M_\beta)$ \cite{BF,Behrend}.

The contraction $\rho:Y\to(\C^2/G)\times\C$ is projective. A map from a connected proper curve to its affine target is constant; for a $T$-fixed stable map, this constant is the origin, the unique fixed point of the target. Thus $\mathcal M_\beta^T$ is the fixed substack of the proper moduli stack of stable maps to $\rho^{-1}(0)$ and is proper.
 For each connected component $F\subset\mathcal M_\beta^T$, the fixed part of the absolute obstruction theory induces $[F]^{\rm vir}$. The moving part of its dual is the virtual normal bundle $N_F^{\rm vir}$, whose equivariant Euler class is invertible after inverting $\nu$. Let $\st_F:F\to\M_{g,n}$ and $\ev_{i,F}:F\to Y$ be the restricted stabilization and evaluation maps. Virtual localization \cite{GP} defines the $T$-equivariant Gromov--Witten CohFT
\begin{equation}\label{eq:cohft}
 \Omega^Y_{g,n}(a_1,\ldots,a_n)
 :=\sum_\beta q^\beta\sum_{F\subset\mathcal M_\beta^T}
 (\st_F)_*\left(
  \frac{\prod_i\ev_{i,F}^*a_i}{e_T(N_F^{\rm vir})}
                       \cap[F]^{\rm vir}\right)\in H^*(\M_{g,n})\otimes\C(\nu)[[q_1,\ldots,q_l]].
\end{equation}

We use the same localization convention for all brackets with superscript $Y$, inserting cotangent-line classes for descendants. The virtual gluing and forgetful identities \cite{Behrend} give the CohFT axioms after localization, with node contraction $\eta^{-1}$ and flat unit $1\in H$. Degree-zero maps give $\Omega^Y_{0,3}(a,b,1)=\eta(a,b)$; the positive-degree terms vanish by the unit axiom.

For a $T$-representation $W_0$ with no zero weights, let $\Omega^{\pt,W_0}$ denote the $T$-equivariant Gromov--Witten theory of a point twisted by the inverse equivariant Euler class \cite[Section~1]{CG}. Its state space is $\C(\nu)$, with unit $1$. On the universal curve $\pi:\mathcal C_{g,n}\to\M_{g,n}$, let $\mathbb E=\pi_*\omega_\pi$ be the rank-$g$ Hodge bundle. We define
\[
 \Omega^{\pt,W_0}_{g,n}(1,\ldots,1)
 :=e_T(R\pi_*\OO\otimes W_0)^{-1}
 =\frac{e_T(\mathbb E^\vee\otimes W_0)}{e_T(W_0)},
\]
where $R\pi_*\OO=\OO-\mathbb E^\vee$ in K-theory. Its pairing is $\eta_{W_0}(1,1)=e_T(W_0)^{-1}$.

Under fixed-locus restriction the state space decomposes as
\[
 H\simeq H_T^*(E)_{\rm loc}
       \oplus\bigoplus_{a\in\{1,\ldots,l-3,L,R\}}\C(\nu)1_{p_a}.
\]
Let $1_E$ restrict to one on $E$ and to zero on every isolated fixed point. The units of the fixed-component blocks sum to the global unit:
\[
 1=1_E+\sum_{a\in\{1,\ldots,l-3,L,R\}}1_{p_a}.
\]

\begin{proposition}\label{prop:split}
The specialization of \eqref{eq:cohft} at $\bq=0$ splits as
\begin{equation}\label{eq:boundary-splitting}
 \Omega^Y\big|_{\bq=0}
 =\Omega^{\PP^1,\,\OO(-2)_{2\nu}\oplus\OO_{-2\nu}}
  \oplus\bigoplus_{a\in\{1,\ldots,l-3,L,R\}}\Omega^{\pt,T_{p_a}Y}.
\end{equation}
The curve summand is twisted by $e_T(R\pi_*f^*V)^{-1}$, and each point summand uses $W_0=T_{p_a}Y$.
\end{proposition}
\begin{proof}

Setting $\bq=0$ retains the degrees $d[E]$. A fixed map lies in the exceptional fiber, and its nonconstant components contribute a nonnegative sum of the independent classes $[E_i]$. For $d>0$, all such components therefore map to $E$; connectedness puts the contracted components there as well.
 Its moving deformation complex is $R\pi_*f^*V$. The two normal characters are nonzero, so its fixed deformation theory is that of maps to $E$. In degree zero the connected image lies in one component of $Y^T$. All classes with inputs in distinct summands vanish. At a node, the normalization exact sequence subtracts the fiber of $V$ from the component complexes. The inverse Euler twist consequently contributes $e_T(V)$, as required by the inverse twisted metric. This gives the stated CohFT decomposition.
\end{proof}

\subsection{The rank-two quantum algebra at \texorpdfstring{$\bq=0$}{q=0}}
Let $h\in H^2(E)$ satisfy $\int_Eh=1$. We also write $h$ for its extension by zero to the other fixed components. Since $h^2=0$ for the ordinary cup product on $E$,

\[
 e_T(V)=(2\nu-2h)(-2\nu),\qquad
 \frac1{e_T(V)}=-\frac1{4\nu^2}\left(1+\frac h\nu\right).
\]
Integrating the appropriate products against this inverse Euler class gives the curve-block metric

\begin{equation}\label{eq:curve-metric}
 \bigl(\eta(1_E,1_E),\eta(1_E,h),\eta(h,h)\bigr)
   =\left(-\frac1{4\nu^3},-\frac1{4\nu^2},0\right).
\end{equation}
Bryan--Gholampour's multiple-cover formula \cite[Lemma~3]{BG}, multiplied by the extra-line factor $(-2\nu)^{-1}$, gives $-1/d^3$ in class $d[E]$.
Each divisor insertion of $h$ contributes $\int_{dE}h=d$. Thus the positive-degree three-divisor series is $-\sum_{d\geq1}Q^d=-Q/(1-Q)$, while the degree-zero three-divisor term vanishes. The unit axiom and the divisor equation give
\[
 \eta(h\star h,1_E)=0,\qquad
 \eta(h\star h,h)=-\frac{Q}{1-Q}.
\]
Writing $h\star h=A1_E+Bh$ and using \eqref{eq:curve-metric} determines $A$ and $B$:

\begin{equation}\label{eq:curve-product}
 h\star h=\frac{4\nu^2Q}{1-Q}\,1_E-\frac{4\nu Q}{1-Q}\,h.
\end{equation}
Writing $Q=a^2$, let $e_\pm$ be the primitive idempotents and choose square roots $v_\pm$ of their self-pairings: $v_\pm^2=\eta(e_\pm,e_\pm)$. Explicitly,
\begin{align}
 e_\pm&=\frac{1\pm a}{2}\,1_E
             \ \pm\frac{1-a^2}{4\nu a}\,h,\label{eq:curve-idempotents}\\
 v_+^2&=-\frac{(1+a)^2}{16\nu^3a},\qquad
 v_-^2=\frac{(1-a)^2}{16\nu^3a}.\label{eq:curve-norms}
\end{align}
The classes $\eh_\pm=e_\pm/v_\pm$ form an orthonormal basis of the curve summand. The eigenvalues of quantum multiplication by $h$ are $2\nu a/(1+a)$ and $-2\nu a/(1-a)$. Each isolated fixed point $p$ contributes a one-dimensional summand $\C(\nu)1_p$ in Proposition~\ref{prop:split}, where $1_p$ restricts to one at $p$ and zero on the other fixed components. For $p$ at distance $m$ along an arm, its self-pairing is
\[
 \eta(1_p,1_p)=\frac{1}{e_T(T_pY)}
             =\frac{1}{8m(m+1)\nu^3}.
\]
The partial-boundary algebra is thus semisimple over $\C(\nu)((a))$. At the full large-radius limit the curve relation becomes $h^2=0$.

\subsection{Specialization of the calibration}

We now pass from the CohFT decomposition in Proposition~\ref{prop:split} to its reconstruction matrix. Lemma~\ref{lem:basechange} shows that this matrix commutes with coefficientwise specialization.

On the semisimple locus let $e_\alpha$ be the primitive idempotents, choose square roots $v_\alpha^2=\eta(e_\alpha,e_\alpha)$, and let $u_\alpha$ be the corresponding canonical coordinates. Our canonical fundamental solution is $\Psi R(z)e^{U/z}$ \cite[Section~1]{GiventalHigher}, where $\Psi$ has columns $\eh_\alpha=e_\alpha/v_\alpha$ and $U=\diag(u_\alpha)$. In this normalized frame, stable-graph reconstruction \cite[Sections~2.2--2.3]{MOPPZ} places $R^{-1}(\psi_i)$ at each marking and the bivector
\[
 \frac{I-R^{-1}(z)\bigl(R^{-1}(w)\bigr)^T}{z+w}
\]
at each edge, where $z,w$ are the cotangent classes of its two branches. The symplectic identity $R(z)R(-z)^T=I$ makes this quotient a power series. The flat-unit translation is $z(v-R^{-1}(z)v)$, where $v=(v_\alpha)_\alpha$ represents the unit.

\begin{lemma}\label{lem:basechange}
Let a flat-unit CohFT over an integral $\mathbb Q$-algebra $A$ have a split semisimple Frobenius algebra with idempotent basis $e_\alpha$. Choose $v_\alpha\in A^\times$ with $v_\alpha^2=\eta(e_\alpha,e_\alpha)$. In the basis $e_\alpha/v_\alpha$, let $R(z)=I+\sum_{k\geq1}R_kz^k$ be its reconstruction matrix, initially taken over $\operatorname{Frac}(A)$. Each $R_k$ has entries in $A$. A specialization $A\to B$ preserving these hypotheses sends $R$ to the reconstruction matrix of the specialized CohFT.
\end{lemma}
\begin{proof}

Over $\operatorname{Frac}(A)$, write $R^{-1}(z)=I+\sum_{k\geq1}C_kz^k$. Fix $k\geq1$ and choose $g>3k$. We use the pointed classes
\[
 \kappa_j=(\pi_{n+1})_*(\psi_{n+1}^{j+1})
\]
on $\mathcal M_{g,n}$, obtained by restricting the usual classes from the compactification, where $\pi_{n+1}:\M_{g,n+1}\to\M_{g,n}$ forgets the last marking. This fixes our convention for the pointed logarithmic $\kappa$-classes.

The stable-range coefficient-recovery argument of \cite[Theorem~2 and Lemma~2.2]{MOPPZ} applies with these generators. On $\mathcal M_{g,2}$ only the single-vertex graph survives. Nontrivial translation terms contribute polynomials of positive degree in the pointed $\kappa$-classes, so the pure $\psi_1^k$ coefficient comes from $C_k$ on the first leg and the identity on the second. It is
\[
 [\psi_1^k]\,
 \Omega_{g,2}(e_\alpha,e_\beta)\big|_{\mathcal M_{g,2}}
 =
 \eta(e_\beta,e_\beta)^{1-g}
 \frac{v_\alpha}{v_\beta}(C_k)^\beta{}_\alpha.
\]
Here $[\psi_1^k]$ extracts the pure monomial $\psi_1^k$. By stable-range independence this is a rational linear functional on the tautological subspace, which we extend rationally to the ambient cohomology. The ratio $v_\alpha/v_\beta$ converts the matrix from the normalized basis to the ordinary idempotent basis.

The class on the left is defined over $A$, and the displayed norm and ratio are units in $A$. Hence every $C_k$ belongs to $A$, and formal inversion gives $R_k\in A$. The reconstruction identity, initially over $\operatorname{Frac}(A)$, therefore holds over $A$. In fixed cohomological degree its graph and translation sums are finite, so the identity commutes with $A\to B$. The specialized matrix reconstructs the specialized CohFT and is its matrix by uniqueness. The matrix called $R$ in \cite{MOPPZ} is the inverse of ours.

\end{proof}

Regrouping Novikov series in powers of $Q=a^2$ gives a coefficientwise map to
\begin{equation}\label{eq:formal-ring}
 A_0=\C(\nu)((a))[[\bq]].
\end{equation}
Set $K'=\C(\sqrt\nu)((\sqrt a))$ and $A=K'[[\bq]]$. This finite extension of the coefficient field contains the square roots of the boundary idempotent norms.
The ring $A$ is complete in the ideal generated by $\bq$. The split boundary idempotents lift uniquely by idempotent lifting, and their norms remain units. Each chosen square root lifts by Hensel's lemma, since its derivative $2v_\alpha$ is a unit. Thus the hypotheses of Lemma~\ref{lem:basechange} hold over $A$.
 Applying Lemma~\ref{lem:basechange} to $A\to K'$, $\bq\mapsto0$, with $Q=a^2$ retained, gives
\begin{equation}\label{eq:a-r-limit}
 R^Y\big|_{\bq=0}
   =R^E\oplus\left(\bigoplus_{m=1}^{l-3}R_m^Y\right)
                    \oplus(R_1^Y)^{\oplus2}.
\end{equation}
Here $R^E$ is the calibration of the curve block, and $R_m^Y$ is the scalar calibration of the twisted point with weights \eqref{eq:point-weights}. \rev{The final two copies correspond to $p_L$ and $p_R$.}

Let $C_{D_l}=((\alpha_i,\alpha_j))$ be the Cartan matrix, and let $D_i=-\sum_j(C_{D_l}^{-1})_{ji}[E_j]$ be the equivariant divisor classes, pulled back to $Y$, satisfying $\int_{E_j}D_i=\delta_{ij}$. The unit and these divisors form a flat basis of $H$. Write a formal primary parameter as
\[
 \tau=t_0\,1+\sum_{i=1}^l t_iD_i.
\]
For $2g-2+n>0$, the primary shift is
\begin{equation}\label{eq:primary-shift}
 \Omega^{Y,\tau}_{g,n}(a_1,\ldots,a_n)
 =\sum_{m\geq0}\frac1{m!}(\pi_m)_*
  \Omega^Y_{g,n+m}(a_1,\ldots,a_n,\underbrace{\tau,\ldots,\tau}_{m}),
\end{equation}
where $\pi_m:\M_{g,n+m}\to\M_{g,n}$ forgets the extra markings. For a class $\beta=\sum_i d_i[E_i]$, the divisor equation changes its weight to
\[
 q^\beta\exp\!\left(\sum_i d_it_i\right)
       =\prod_i(q_i e^{t_i})^{d_i}.
\]
Hence the divisor part of the shift is absorbed by $q_i\mapsto q_i e^{t_i}$. For the shifted theory, the mirror parameters are obtained by this substitution.

Let $\star_\tau$ be the equivariant quantum product of the shifted theory. Its quantum multiplication one-form is the $\End(H)$-valued form
\[
 C=I\,dt_0+\sum_{i=1}^l C_i(\tau)\,dt_i,
 \qquad C_i(\tau)(a)=D_i\star_\tau a.
\]
Teleman's deformation equations \cite[Propositions~7.7 and 8.2]{Teleman} show that the calibration selected by the shifted CohFT satisfies
\begin{equation}\label{eq:qde}
 z\,d(\Psi R^Y e^{U/z})=C\,\Psi R^Y e^{U/z},
\end{equation}
where $d$ differentiates the primary coordinates with the Novikov symbols fixed. For $F=\Psi R^Y e^{U/z}$, its components are $z\partial_{t_i}F=(D_i\star_\tau)F$ and $z\partial_{t_0}F=F$.

\section{The Toda curve and its calibration}\label{sec:mirror}

We describe the local mirror data and construct its stationary-phase matrix. The genus-zero mirror theorem identifies the Frobenius structures; Proposition~\ref{prop:b-flat} places the B-model matrix in the same quantum differential system as the A-model calibration.

\subsection{The mirror map and metric}
Earlier relativistic Toda constructions beyond type $A$, including an affine type-$D$ example, appear in Kruglinskaya--Marshakov \cite{KrM}. Related Landau--Ginzburg models for extended affine Weyl Frobenius manifolds were constructed by Dubrovin--Strachan--Zhang--Zuo \cite{DSZZ} and Brini--van Gemst \cite{BvG}.

Throughout this paper, ``Toda curve'' denotes the degeneration of the type-$D_l$ affine relativistic Toda spectral curve at vanishing base Novikov parameter. Its rational normalization is $\PP^1_\mu$, and its equation is
\begin{equation}\label{eq:toda}
 \lambda(\mu)=K
 \frac{\prod_{i=1}^l(\mu-\kappa_i)(\mu-\kappa_i^{-1})}
      {\mu^{l-2}(\mu^2-1)^2}
 =K\frac{P(w)}{w^2-4},\quad
 w=\mu+\mu^{-1},\quad P(w)=\prod_{i=1}^l(w-c_i),
\end{equation}
where $c_i=\kappa_i+\kappa_i^{-1}$.
The partial large-radius boundary used below for calibration is a further specialization within the K\"ahler moduli of $Y$, obtained by setting $q_i=0$ for $i\ne l-2$ while retaining $Q=q_{l-2}\ne0$.
Put
\begin{equation}\label{eq:phase}
 W=2\nu\log\lambda,\qquad \xh=-W,\qquad
 \phi_D=\sqrt\nu\,d\log\mu,\qquad \iota(\mu)=\mu^{-1}.
\end{equation}
The mirror map is
\begin{equation}\label{eq:mirror-map}
 q_i=\frac{\kappa_{i+1}}{\kappa_i}\ (i<l),\qquad
 q_l=\frac1{\kappa_{l-1}\kappa_l},\qquad t_0=2\nu\log K.
\end{equation}
We first work where $P$ has simple roots different from $\pm2$ and
\begin{equation}\label{eq:critical-polynomial}
 N(w)=(w^2-4)P'(w)-2wP(w)
\end{equation}
has $l+1$ simple roots with distinct critical values of $\lambda$. Each root $w_\alpha$ has two reciprocal lifts under $\mu\mapsto w$; choose one and denote it by $p_\alpha$. On a simply connected domain choose logarithms of $\lambda$ and Morse coordinates on the quotient:
\begin{equation}\label{eq:morse}
 \xh=-u_\alpha+s_\alpha^2,\qquad u_\alpha=W(p_\alpha).
\end{equation}
With these logarithms, $u_\alpha$ is independent of the chosen lift.

The genus-zero mirror theorem for $Z$ is \cite[Theorem~3.2]{BMS}. The surface metric and cubic tensor acquire the common factor $(-2\nu)^{-1}$ on adjoining the extra line. For tangent vectors $X,Y$, put
\[
 \eta_0(X,Y)=\sum_{\text{critical lifts }p}
 \Res_p\frac{(X\log\lambda)(Y\log\lambda)\phi_D^2}{d\log\lambda}.
\]
The threefold metric is $-\eta_0/(4\nu^2)$: contraction of the surface cubic with the unit $(2\nu)^{-1}\partial_{\log K}$ first contributes $(2\nu)^{-1}$. Summing over the two lifts of a critical orbit gives
\begin{equation}\label{eq:metric-local}
 v_\alpha^2=\eta(e_\alpha,e_\alpha)
   =\frac1{8\nu^2}
       \left(\left.\frac{d\log\mu}{ds_\alpha}\right|_{p_\alpha}\right)^2.
\end{equation}
Indeed $e_\alpha\log\lambda=1/(2\nu)$ at the corresponding critical orbit. Each lift contributes
\[
 -\frac14\left(\frac{d\log\mu}{ds_\alpha}\right)^2
 \quad\text{to }\eta_0(e_\alpha,e_\alpha).
\]

\subsection{Equivariant recursion and stationary phase}\label{subsec:equivariant-recursion}
Let $\mathcal L_{\rm sgn}$ be the flat sign line bundle on $Q^\circ=\PP^1_w\setminus\{\pm2\}$ associated to the reciprocal double cover. It has monodromy $-1$ around each of $\pm2$ and a canonical pairing $\mathcal L_{\rm sgn}^{\otimes2}\simeq\OO_{Q^\circ}$. Choosing $p_\alpha$ fixes a local flat frame; exchanging it with its inverse reverses that frame. The odd differential $\phi_D$ descends as a meromorphic $\mathcal L_{\rm sgn}$-valued one-form. On the chosen critical disks use
\begin{equation}\label{eq:y-kernel}
 \yt_\alpha=\frac{i}{2\nu}
       (\log\mu-\log\mu(p_\alpha)),\qquad
 B^-=B^{\rm std}-(\id\times\iota)^*B^{\rm std},\qquad
 B^{\rm std}=\frac{d\mu_1d\mu_2}{(\mu_1-\mu_2)^2}.
\end{equation}
Here $\yt_\alpha$ is odd under exchange of lifts, so it is a local section of $\mathcal L_{\rm sgn}$; $B^-$ is odd in each variable and is an $\mathcal L_{\rm sgn}\boxtimes\mathcal L_{\rm sgn}$-valued bidifferential. On all $2(l+1)$ lifted disks, $(\xh,\yt,B^-/2)$ are the local $\mathbb Z_2$-equivariant data of Giacchetto--Kramer--Lewa\'nski \cite[Section~8]{GKL}: $d\xh$ has trivial character, while $\yt$ and the kernel have sign character. We reduce to one disk per orbit and use the unhalved $B^-$, with standard diagonal principal part, as in \cite{FLM}. The following formulas express this line-bundle-valued recursion in the chosen frames. The local conjugation $p\mapsto\bar p$ is $s_\alpha\mapsto-s_\alpha$.

Start with $\omega_{0,1}=\yt\,d\xh$ and $\omega_{0,2}=B^-$. For $2g-2+n>0$, define
\begin{align}\label{eq:recursion}
 \omega_{g,n}(p_1,\mathbf p)
 =\sum_\alpha\Res_{p=p_\alpha}K_\alpha(p_1,p)
 \Bigl[&\omega_{g-1,n+1}(p,\bar p,\mathbf p)\\[-2pt]
 &+\sum_{g_1+g_2=g,\ I\sqcup J=\{2,\ldots,n\}}^{\prime}
 \omega_{g_1,|I|+1}(p,p_I)\omega_{g_2,|J|+1}(\bar p,p_J)\Bigr],\notag
\end{align}
where the prime excludes factors $\omega_{0,1}$, and
\begin{equation}\label{eq:recursion-kernel}
 K_\alpha(p_1,p)=
 \frac{\tfrac12\int_{\bar p}^{p}B^-(p_1,\cdot)}
      {(\yt_\alpha(p)-\yt_\alpha(\bar p))d\xh(p)}.
\end{equation}
The genus $g-1$ term is omitted for $g=0$. The stable forms $\omega_{g,n}$ are meromorphic sections of $(K_{Q^\circ}\otimes\mathcal L_{\rm sgn})^{\boxtimes n}$. If $\omega^{G,\rm av}_{g,n}$ denotes GKL's recursion over all lifted critical points with kernel $B^-/2$, then $\omega_{g,n}=2^{g+n-1}\omega^{G,\rm av}_{g,n}$ on the representative disks, by \cite[Proposition~8.8 and Remark~2.2]{GKL}.

For the calibration put $\rho_\alpha=\sqrt2s_\alpha$ and expand
\begin{equation}\label{eq:evaluated-kernel}
 \mathcal B_\beta(p)=
  \left.\frac{B^-(p,p')}{d\rho_\beta(p')}\right|_{p'=p_\beta}
  =\left(\frac{\delta_{\alpha\beta}}{\rho_\alpha^2}
            +\sum_{j\geq0}b_j^{\beta\alpha}\rho_\alpha^j\right)d\rho_\alpha.
\end{equation}
Define
\begin{equation}\label{eq:r-jets}
 R^B_{\beta\alpha}(z)=\delta_{\alpha\beta}
       -\sum_{j\geq0}(2j-1)!!\,b_{2j}^{\beta\alpha}z^{j+1}.
\end{equation}
To express this series by stationary phase, choose a ray for $z$ and a compatible square root $\sqrt z$. Let $\Gamma_\alpha$ denote the local steepest-descent contour through $p_\alpha$, oriented so that $\rho_\alpha/\sqrt z$ increases along the real axis. The integral below is interpreted by its formal Gaussian expansion; $\int^{\rm reg}$ detours around the possible double pole at $\rho_\alpha=0$. With this convention, \eqref{eq:r-jets} is the expansion of

\begin{equation}\label{eq:r-integral}
 -\frac{\sqrt z}{\sqrt{2\pi}}
 \int_{\Gamma_\alpha}^{\rm reg}
       e^{-\rho_\alpha^2/(2z)}\mathcal B_\beta.
\end{equation}
Integration of $d(e^{-\rho^2/(2z)}/\rho)$ fixes the double-pole contribution to be $1$. Each coefficient of \eqref{eq:r-jets} depends on finitely many derivatives of $\xh$ and of the kernel. It therefore extends wherever the critical points remain simple and avoid the zeros and poles. Section~\ref{sec:calibration} realizes the scalar boundary expansions by convergent beta integrals. Changing lift frames by a diagonal sign matrix $D_{\rm sgn}$ sends $R^B\mapsto D_{\rm sgn}R^BD_{\rm sgn}$; the compatible normalized idempotents change by the same signs.

\begin{proposition}\label{prop:b-flat}
The series $R^B$ is symplectic, and $\Psi R^B e^{U/z}$ satisfies \eqref{eq:qde} under \eqref{eq:mirror-map}.
\end{proposition}
\begin{proof}
The Prym double-Laplace factorization \rev{\cite[Proposition~4.4]{FLM}} gives $R^B(z)^TR^B(-z)=I$. For flatness we use the smooth-family calculation of \rev{\cite[Proposition~6.2]{FLM}}, recording its normalization in our notation. For $\tau_i=\log\kappa_i$, the primary differentials are the kernel periods
\[
 \omega_{\tau_i}=-\nu\sqrt\nu\int_{\kappa_i^{-1}}^{\kappa_i}B^-,
 \qquad
 \omega_1=-\frac{\sqrt\nu}{2}\int_0^\infty B^-.
\]
Direct integration verifies these formulas. For a flat basis $T_a$, put $\omega_a=(T_aW)\phi_D$. With compatible square roots, \eqref{eq:metric-local} gives
\[
 (\eta\Psi)_{a\alpha}=\eta(T_a,\eh_\alpha)
 =\frac{1}{2\nu\sqrt\nu}
   \left.\frac{\omega_a}{d\rho_\alpha}\right|_{p_\alpha}.
\]
In the convention $u_\alpha=-\xh(p_\alpha)$, evaluation of the logarithmic Rauch formula gives the rotation coefficients
\[
 \Gamma_{\alpha\beta}
 =\frac{\partial_{u_\alpha}v_\beta}{v_\alpha}
 =-\left.
 \frac{B^-(p,q)}{d\rho_\alpha(p)\,d\rho_\beta(q)}
 \right|_{p=p_\alpha,\,q=p_\beta}
 \quad(\alpha\ne\beta),\qquad
 \Gamma_{\alpha\alpha}=0.
\]
The minus sign comes from replacing the critical values of $\xh$ by $u_\alpha$. Applying the same variation to the primary differentials gives
\[
 d(\eta\Psi)=-\eta\Psi[\Gamma,dU],
 \qquad d\Psi=-\Psi[\Gamma,dU],
\]
where the second identity uses the constancy of $\eta$.
The relative paths acquire no endpoint term: in a moving coordinate $\xi$ at a logarithmic puncture, $d\xh=c\,d\xi/\xi+O(d\xi)$, and the horizontal velocity at fixed $\xh$ is $O(\xi)$. The kernel endpoint evaluation therefore vanishes.
Differentiating \eqref{eq:r-integral} and integrating by parts now gives

\[
 dR^B=z^{-1}[dU,R^B]+[\Gamma,dU]R^B.
\]
The genus-zero mirror identification gives $C\Psi=\Psi\,dU$. Hence
\[
 \begin{aligned}
 z\,d(\Psi R^Be^{U/z})
 &=\Psi[dU,R^B]e^{U/z}
       +\Psi R^B\,dU\,e^{U/z}\\
 &=C\,\Psi R^Be^{U/z},
 \end{aligned}
\]
which proves the assertion.

\end{proof}

\subsection{The partial large-radius degeneration}\label{sec:degeneration}
We now study the mirror data at the partial boundary used in Section~\ref{sec:amodel}. The different scales of the zeros separate the critical points into one-dimensional blocks and a two-dimensional block. We first describe the limiting curve, then take the limit of its local kernels and stationary-phase matrix.

\subsubsection{The semistable curve}
Fix $Q\neq0,1$. An approach to the partial boundary is given by
\begin{equation}\label{eq:arc}
 \begin{split}
 \kappa_i&=A_i t^{-d_i}\quad(1\leq i\leq l-3),\\
 \kappa_{l-2}&=A t^{-d_f},\qquad
 \kappa_{l-1}=AQ t^{-d_f},\qquad
 \kappa_l=Ct^{-d_l},
 \end{split}
 \qquad d_1>\cdots>d_{l-3}>d_f>d_l>0.
\end{equation}
The constants are nonzero and the exponents can be taken integral after base change. By \eqref{eq:mirror-map}, $q_{l-2}=Q$ and every other $q_i$ tends to zero.

\begin{proposition}\label{prop:semistable}
Mark the zero and pole sections of $\lambda$ on the quotient $\PP^1_w$. The one-parameter family defined by \eqref{eq:arc} admits a semistable model extending across $t=0$. After contraction of unmarked two-point components, the central fiber (the curve over $t=0$) is a union of rational curves arranged in the chain
\begin{equation}\label{eq:component-chain}
 C_0\text{---}C_l\text{---}C_f\text{---}C_{l-3}\text{---}\cdots\text{---}C_1.
\end{equation}
Each dash represents one node joining adjacent components. On each component $C_j$, define the meromorphic function $\Lambda_j=\lim_{t\to0}\lambda/A_j(t)$ in the coordinate shown below, where the nonzero scalar factors $A_j(t)$ are specified in the proof. The last column lists the zeros of $d\log\Lambda_j$ in the coordinate given in the second column.
\begin{equation}\label{eq:component-table}
\begin{array}{c|c|c|c}
 \text{component}&\text{coordinate}&\Lambda_j&\text{critical points}\\ \hline
 C_0&w&(w^2-4)^{-1}&0\\
 C_l&s=w/c_l&(s-1)/s^2&2\\
 C_f&s=w/c_{l-2}&(s-1)(s-Q)/s&\pm\sqrt Q\\
 C_k,\ 1\leq k\leq l-3&s=w/c_k&s^{m_k}(s-1)&m_k/(m_k+1)
\end{array}
\quad m_k=l-k-2.
\end{equation}
At each attachment, infinity on the smaller-scale component meets zero on the next component. The zero sections lie at $s=1$ on the ordinary tails and at $s=1,Q$ on $C_f$. The poles $w=\pm2$ lie on $C_0$, and the infinity pole lies on $C_1$.

Upstairs, the central component is $\widetilde C_0=\PP^1_\mu$, with $w=\mu+\mu^{-1}$. Every tail in \eqref{eq:component-chain} has two copies. One chain attaches to $\mu=0$ and the other to $\mu=\infty$; the involution exchanges them.
\end{proposition}
\begin{proof}
The scale hierarchy is $1\ll c_l\ll c_{l-2}\asymp c_{l-1}\ll c_{l-3}\ll\cdots\ll c_1$, with $c_{l-1}/c_{l-2}\to Q$. Successively blow up the collisions of the zero sections at infinity. The functions $w$, $w/c_l$, $w/c_{l-2}$, and $w/c_k$ become coordinates on $C_0$, $C_l$, $C_f$, and $C_k$ in \eqref{eq:component-chain}, respectively. For neighboring scales $b\ll c$, coordinates $u=b/w$ and $v=w/c$ satisfy $uv=b/c$. Up to an invertible factor this is $uv=t^d$, where $d>0$ is their exponent difference. Resolving it inserts $d-1$ two-point rational components and produces a regular surface with reduced nodal central fiber. Contracting those components recovers the stated chain. Each displayed component has at least three nodes or marked points.

The reciprocal double cover is branched only at $w=\pm2$, which remain on $C_0$. Its restriction to every other component splits into two copies. The two points over $w=\infty$ on $\widetilde C_0$ are $\mu=0$ and $\infty$, proving the upstairs description.

For the component functions divide $\lambda$ on each chart by, respectively,
\begin{align*}
 A_0(t)&=K(-1)^l\prod_i c_i,&
 A_l(t)&=K(-1)^{l-1}\frac{\prod_{i<l}c_i}{c_l},\\
 A_f(t)&=K(-1)^{l-3}\left(\prod_{i\leq l-3}c_i\right)c_{l-2},&
 A_k(t)&=K(-1)^{k-1}\left(\prod_{i<k}c_i\right)c_k^{l-k-1}.
\end{align*}
For example,
\[
 \frac{\lambda(c_ks)}{A_k(t)}
 =s^{l-k-2}(s-1)
 \frac{\prod_{i<k}(1-c_ks/c_i)\prod_{i>k}(1-c_i/(c_ks))}
      {1-4/(c_k^2s^2)}.
\]
The fraction on the right tends to $1$ on compact subsets of $\PP^1_s\setminus\{0,1,\infty\}$. The same factorization on the other charts gives the table. These limits describe the logarithmic differential on the nodal source: $d\log\lambda\to d\log\Lambda_j$ on each component away from the zeros and poles of $\Lambda_j$. On each critical disk, subtracting the critical value of $\xh$ gives
\[
 \xh(p)-\xh(p_\alpha)
 =-2\nu\log\frac{\lambda(p)}{\lambda(p_\alpha)}.
\]
The scalar factor $A_j(t)$ cancels in this ratio. At each node the two logarithmic residues sum to zero.
\end{proof}

\subsubsection{Critical points and the limiting kernel}
All the critical points in \eqref{eq:component-table} are simple and lie away from the attaching nodes and the zero and pole sections. Uniform convergence of the logarithmic derivatives on small disks about these points, followed by Rouch\'e's theorem, gives $1+1+2+(l-3)=l+1$ nearby critical points. The polynomial \eqref{eq:critical-polynomial} has degree $l+1$ and leading coefficient $l-2$, so these exhaust the critical points.

On a tail $w=cs$ the two lifts have expansions
\[
 \mu_+=cs(1+O(c^{-2})),\qquad
 \mu_-=c^{-1}s^{-1}(1+O(c^{-2})).
\]
Thus on the positive lift
\[
 d\yt\longrightarrow \frac{i}{2\nu}\frac{ds}{s},\qquad
 B^-\longrightarrow \frac{ds_1ds_2}{(s_1-s_2)^2}.
\]
The signs change according to the chosen lifts. For different scales $b\ll c$, the kernel is $O(b/c)+O((bc)^{-1})$, and it is $O(c^{-1})$ between the central chart and a tail. These estimates hold on compact critical disks; Cauchy's formula gives convergence of every fixed-order jet. On $\widetilde C_0$ the kernel remains $B^-$.

\begin{corollary}\label{cor:b-r-limit}
Coefficientwise in $z$, the stationary-phase matrix has limit
\begin{equation}\label{eq:b-r-limit}
 R^B\longrightarrow R^B_0\oplus R^B_l\oplus R^B_f
                          \oplus\bigoplus_{k=1}^{l-3}R^B_k.
\end{equation}
The central and ordinary-tail blocks have rank one; $R^B_f$ has rank two.
\end{corollary}
\begin{proof}
The difference $\xh(p)-\xh(p_\alpha)$ is independent of $A_j(t)$. The nonzero limiting Hessian implies convergence of the Morse coordinate and its inverse. Now apply the finite-jet formula \eqref{eq:r-jets} to the preceding kernel limits.
\end{proof}

\subsubsection{The full large-radius limit}
The further limit $Q=a^2\to0$ explains the exceptional role of the partial boundary. On the component $C_f$,
\[
 \Lambda_f(a\zeta;a^2)=-1+a(\zeta+\zeta^{-1})-a^2.
\]
Separating the zero sections at $s=1,a^2$ and the critical sections at $s=\pm a$ produces three components with coordinates $v=s/a^2$, $\zeta=s/a$, and $s$. Their functions tend to $v^{-1}-1$, $-1$, and $s-1$, respectively. The two critical sections meet the middle component at $\zeta=\pm1$. On this component the limiting differential is zero:
\[
 d\xh\longrightarrow0,\qquad
 \frac{d\xh}{2\nu a}\longrightarrow d(\zeta+\zeta^{-1}).
\]
This collision corresponds to the relation $h^2=0$ in \eqref{eq:curve-product}. The comparison below is made at $Q\neq0$, where both critical points remain Morse.

\subsection{The boundary calibration}\label{sec:calibration}
We identify the blocks of Corollary~\ref{cor:b-r-limit} with the A-model blocks in \eqref{eq:a-r-limit}. Beta integrals compute the scalar blocks, while the toric calibration theorem identifies the rank-two block. \rev{These boundary values will fix the constant diagonal ambiguity in the quantum differential system.}

\subsubsection{The point blocks}
For $m\geq1$, inverse-Euler quantum Riemann--Roch \cite[Theorem~1 and Section~4]{CG} and \eqref{eq:point-weights} give
\begin{equation}\label{eq:point-r}
 R_m^Y(z)=\exp\left[
 \sum_{j\geq1}\frac{B_{2j}}{2j(2j-1)}
 \{m^{1-2j}+1-(m+1)^{1-2j}\}
 \left(\frac z{2\nu}\right)^{2j-1}\right].
\end{equation}
Indeed the inverse-Euler point formula has exponent $-\sum_{j\geq1}B_{2j}z^{2j-1} \sum_b w_b^{1-2j}/[2j(2j-1)]$ for tangent weights $w_b$.

For an ordinary tail $C_k$, put $m=m_k=l-k-2$. Its function in \eqref{eq:component-table} is $\Lambda_k(s)=s^m(s-1)$. Its logarithmic critical point is $p=m/(m+1)$ and $\xh''(p)=2\nu(m+1)^3/m$. Set $X=2\nu/z$. Formula \eqref{eq:metric-local} gives its norm as $1/[4\nu^2p^2\xh''(p)]=1/[8m(m+1)\nu^3]$, agreeing with the point theory. The evaluated kernel is $ds/[\sqrt{\xh''(p)}(s-p)^2]$. Since
\[
 \frac{\xh'(s)}{s-p}=\frac{2\nu(m+1)}{s(1-s)},
\]
put
\[
 E_X(s)=e^{-(\xh(s)-\xh(p))/z}
       =\left(\frac{s^m(1-s)}{p^m(1-p)}\right)^X.
\]
The identity
\[
 d\left(\frac{E_X(s)}{s-p}\right)
 =-\frac{\xh'(s)E_X(s)}{z(s-p)}\,ds
   -\frac{E_X(s)}{(s-p)^2}\,ds
\]
removes the double pole from \eqref{eq:r-integral}. On $(0,1)$, with $\Re X>0$, the endpoint terms vanish and the regularized integration by parts gives the beta expression
\[
 \sqrt{\frac{mX}{2\pi(m+1)}}
 \left(\frac{(m+1)^{m+1}}{m^m}\right)^X
 \int_0^1s^{mX-1}(1-s)^{X-1}\,ds.
\]
Using $\mathrm B(mX,X)=(m+1)\mathrm B(mX,X+1)$ yields

\begin{equation}\label{eq:beta-r}
 R_m^B(z)\sim
 \sqrt{\frac{m(m+1)X}{2\pi}}
 \left(\frac{(m+1)^{m+1}}{m^m}\right)^X
 \frac{\Gamma(mX)\Gamma(X+1)}{\Gamma((m+1)X+1)}.
\end{equation}
Here $\sim$ denotes the stationary-phase expansion as $z\to0$. In Stirling's formula, the leading exponential and square-root factors of the Gamma quotient cancel the prefactors in \eqref{eq:beta-r}. The remaining logarithm is
\[
 \log R_m^B(z)
 =\sum_{j\geq1}\frac{B_{2j}}{2j(2j-1)}
 \bigl(m^{1-2j}+1-(m+1)^{1-2j}\bigr)X^{1-2j}.
\]
Since $X=2\nu/z$, this is the exponent in \eqref{eq:point-r}, and $R_m^B=R_m^Y$.

On $C_0$, the substitution $u=\mu^2/(\mu^2-1)$ turns $(w^2-4)^{-1}$ into the $m=1$ function $u(u-1)$. At the chosen critical lift $\mu=i$, the evaluated Cauchy kernel in $u$ equals the evaluated Prym kernel, since
\[
 \frac{du_1du_2}{(u_1-u_2)^2}
 =\left(\frac1{(\mu_1-\mu_2)^2}
        -\frac1{(\mu_1+\mu_2)^2}\right)d\mu_1d\mu_2,
 \qquad
 \frac1{(i\mu-1)^2}=-\frac1{(\mu+i)^2}.
\]
On $C_l$, the M\"obius substitution $u=1/s$ gives $u(1-u)$. The factor $-1$ cancels after dividing the function by its value at the critical point. Thus the ordinary tails contribute one block for each $m=1,\ldots,l-3$, and $C_0$ and $C_l$ each contribute an additional $m=1$ block. These are the $l-1$ scalar blocks in \eqref{eq:a-r-limit}.

\subsubsection{The rank-two block}
\rev{The component $C_f$ retains two coupled critical points. We identify its threefold calibration using toric remodeling, after matching the flat parameter and specializing an auxiliary torus action.}

The remaining component is the mirror of $\operatorname{Tot}_{\PP^1}(\OO(-2)\oplus\OO)$:
\begin{equation}\label{eq:fork-curve}
 Y_0+s+\frac Qs-(1+Q)=0,\qquad
 Y_0=-\Lambda_f(s),\qquad \Lambda_f(s)=\frac{(s-1)(s-Q)}s.
\end{equation}
Here $(s,Y_0)$ are coordinates on the ambient torus $(\C^*)^2$. The complex parameter in the normalized equation below is $q=Q/(1+Q)^2$. Near $Q=0$, the limit $q\to0$ is the large complex structure limit, corresponding to the A-model large-radius limit $Q\to0$ \cite[Section~3.9.5]{FLZ}. To compare the B-model calibration with $R^E$, we identify its flat parameter with the $Q$ used in the curve summand of Proposition~\ref{prop:split}.

Put $\xh_f=-2\nu\log\Lambda_f$. Since $\lambda/A_f(t)\to\Lambda_f$, this is the limit of the global $\xh$ after subtracting $-2\nu\log A_f(t)$, which is constant in $s$. In particular, $d_s\xh\to d_s\xh_f$. Since $Y_0=-\Lambda_f$, compatible logarithm branches give $2\nu\log Y_0=-\xh_f+2\nu\log(-1)$.

\begin{lemma}\label{lem:fork-unit}
For the normal weights $(2\nu,-2\nu)$, the mirror map of \eqref{eq:fork-curve}, with superpotential $2\nu\log Y_0$, is $\tau=(\log Q)h$: the flat divisor coordinate is $\log Q$, and the unit coordinate is zero, up to a constant depending on the logarithm branch.
\end{lemma}
\begin{proof}
We use the equivariant toric mirror theorem in the form recalled in \cite[Theorem~3.11]{FLZ}, with its general toric source \cite[Theorem~31]{CCIT}, to read the flat coordinates from the $z^{-1}$ coefficient of the normalized $I$-function. Specializing to the local curve gives
\[
 I(q,z)=e^{h\log q/z}
 \sum_{d\geq0}q^d
 \frac{\prod_{j=0}^{2d-1}(2\nu-2h-jz)}
      {\prod_{j=1}^{d}(h+jz)^2}.
\]
Its coefficient of $z^{-1}$ is $h\log q-(2\nu-2h)f(q)$, where
\[
 f(q)=\sum_{d\geq1}\frac{(2d-1)!}{(d!)^2}q^d
     =-\log\frac{1+\sqrt{1-4q}}2.
\]
On the branch at $Q=0$, $f(q)=\log(1+Q)$, giving $h\log Q-2\nu\log(1+Q)\,1$. The normalized toric equation $\bar Y+\bar s+q/\bar s-1=0$ is obtained from \eqref{eq:fork-curve} by dividing $Y_0,s$ by $1+Q$. Returning to $Y_0$ adds $2\nu\log(1+Q)$ to the superpotential, canceling this unit term.
\end{proof}

Introduce an auxiliary torus weight $\varepsilon$ on the base, and denote the resulting A-model calibration by $R^{E,\varepsilon}$; at $\varepsilon=0$ the theory is the curve summand of Proposition~\ref{prop:split}. The tangent weights at the two fixed points are
\[
 (\varepsilon,2\nu,-2\nu-\varepsilon),\qquad
 (-\varepsilon,2\nu+2\varepsilon,-2\nu-\varepsilon).
\]
The equivariant superpotential on the curve \cite[Section~4.6]{FLZ} is $P_\varepsilon=-\varepsilon\log s+(2\nu+\varepsilon)\log Y_0$. At $\varepsilon=0$ it equals $-\xh_f+2\nu\log(-1)$. Let $\check R_\varepsilon$ be the B-model matrix in the conventions of \cite{FLZ}. Fang--Liu--Zong's calibrated matrix theorem \cite[Theorem~7.1]{FLZ} gives $R^{E,\varepsilon}(z)=\check R_\varepsilon(-z)$ for generic $\varepsilon$. Subtracting $P_\varepsilon(p)$ from $P_\varepsilon$ and replacing $z$ by $-z$ gives the exponential $e^{(P_\varepsilon-P_\varepsilon(p))/z}$. At $\varepsilon=0$ this is $e^{-(\xh_f-\xh_f(p))/z}$. The compactified mirror has the Cauchy kernel.

We now specialize to $\varepsilon=0$. On the A-side the proper base is $\PP^1$ and the inverse-Euler normal characters remain $2\nu,-2\nu$, so the twisted CohFT is regular there. Its product remains semisimple at $Q\neq0,1$, and Lemma~\ref{lem:basechange} gives
\[
 R^{E,\varepsilon}\big|_{\varepsilon=0}=R^E.
\]
On the B-side, the limiting phase is $\xh_f$ and the kernel is Cauchy; the nonzero limiting Hessians allow specialization in the finite-jet formula. Thus $\check R_\varepsilon(-z)\big|_{\varepsilon=0}=R_f^B(z)$, and we obtain

\begin{equation}\label{eq:fork-r-match}
 R_f^B(z)=R^E(z).
\end{equation}
The critical points $s=a,-a$ correspond to $e_-,e_+$ respectively: their critical values, up to a common constant, are $4\nu\log(1-a)$ and $4\nu\log(1+a)$, whose $Q\partial_Q$ derivatives are the corresponding eigenvalues of the quantum multiplication operator $h\star(-)$.

The rank-two matrix is explicit to every order by \eqref{eq:r-jets} with $\Lambda_f$ and the Cauchy kernel. For example, choosing
\[
 \sqrt{\xh_f''(a)}=\frac{2\sqrt\nu}{\sqrt a(1-a)},\qquad
 \sqrt{\xh_f''(-a)}=\frac{2i\sqrt\nu}{\sqrt a(1+a)}
\]
and ordering the critical points as $(a,-a)$ gives
\begin{equation}\label{eq:fork-r1}
 R_f^B(z)=I+\frac{z}{32\nu a}
 \begin{pmatrix}
  1+a^2&2i(1-a^2)\\
  2i(1-a^2)&-1-a^2
 \end{pmatrix}+O(z^2).
\end{equation}
Expansion in the Morse coordinate gives the diagonal entry $(\xh''\xh''''-(\xh''')^2)/(24(\xh'')^3)$; its off-diagonal entry is $-1/[\sqrt{\xh''(a)}\sqrt{\xh''(-a)}(2a)^2]$. The pole at $a=0$ records the loss of semisimplicity at full large radius.

\subsubsection{Equality away from the boundary}
The preceding comparisons identify every boundary block. The standard calibration ambiguity is a constant diagonal odd exponential \cite[Proposition~6.7(b)]{GiventalQuant}. We give the argument over our formal coefficient field and use the divisor equation to make this constancy apply to the Novikov degeneration; the boundary values then determine the remaining series.

\begin{proposition}\label{prop:r-equality}
The matrices $R^Y$ and $R^B$ agree under \eqref{eq:mirror-map}.
\end{proposition}
\begin{proof}
Work over a common finite algebraic extension of $F=\operatorname{Frac}(\C(\nu)((a))[[\bq]])$ containing the critical roots, their lifts, and the metric and Hessian square roots. Every coefficient of $R^B$ belongs to this field. In fact, the mirror parameters are algebraic in $q_i$, and inversion of the Morse coordinate uses only rational operations and the chosen Hessian root. The Taylor series of $\xh-\xh(p_\alpha)$ has algebraic coefficients. The differential $dU$ is algebraic as well.

Set $D=(R^Y)^{-1}R^B$. Equations \eqref{eq:qde} and Proposition~\ref{prop:b-flat} imply
\[
 z\,dD=[dU,D].
\]
Writing $D=I+\sum_{j\geq1}D_jz^j$, the first equation makes $D_1$ diagonal. At the next order the diagonal equation gives $dD_1=0$ and the off-diagonal equation makes $D_2$ diagonal. Induction shows that every $D_j$ is diagonal and constant in all Frobenius directions. Symplecticity further makes $\log D$ odd in $z$.

To use this constancy along a Novikov degeneration, we apply the divisor equation from Section~\ref{sec:amodel}. For $\tau_{\rm div}=\sum_i t_iD_i$, put $qe^t=(q_ie^{t_i})_{i=1}^l$. With compatible normalized idempotents, the shifted CohFT and uniqueness of its calibration give
\[
 R^Y(q,\tau_{\rm div};z)=R^Y(qe^t,0;z).
\]
The same identity holds for $R^B$, since its mirror parameters are obtained by the substitution $q_i\mapsto q_ie^{t_i}$. Consequently
\[
 q_i\partial_{q_i}D(q,0;z)
 =\left.\partial_{t_i}D(q,\tau_{\rm div};z)\right|_{t=0}=0.
\]
Thus $D$ is also constant along the Novikov paths used below.

To evaluate $D$, keep $Q=a^2$ and take $q_j=b_jt^{d_j}$ for $j\neq l-2$, with positive integer exponents satisfying $d_l>d_{l-1}$ and with the $b_j$ algebraically independent over $\C(\nu)((a))$. The mirror map gives
\[
 \kappa_l^2=\frac{q_{l-1}}{q_l},\qquad
 \kappa_{l-1}^2=\frac1{q_{l-1}q_l},\qquad
 \kappa_i=\frac{\kappa_{i+1}}{q_i}\quad(i\leq l-2).
\]
After choosing square roots, these identities give the scale hierarchy of \eqref{eq:arc}. Under the Novikov substitution,

\[
 \prod_jq_j^{n_j}\longmapsto
 \left(\prod_jb_j^{n_j}\right)t^{\sum_jd_jn_j}.
\]
For each fixed power of $t$, only finitely many exponent vectors occur, and the corresponding distinct monomials in the $b_j$ are linearly independent over $\C(\nu)((a))$. Thus the map into $\C(\nu)((a))(\mathbf b)[[t]]$ is injective. After a finite extension $L$ of $\C(\nu)((a))(\mathbf b)$ and a finite base change, this embedding extends to the chosen algebraic field inside $L((t^{1/N}))$.

 The A-model limit is \eqref{eq:a-r-limit}. The finite-jet computations in Subsection~\ref{sec:degeneration} give the B-model limit in the same formal field. Equations \eqref{eq:point-r}--\eqref{eq:fork-r-match} identify the two limits. Reducing $R^B=R^YD$ along this arc gives $R_0=R_0D$. Since $R_0=I+O(z)$ is invertible, $D=I$.
\end{proof}

\section{Descendants and free energies}\label{sec:correlators}
Write $R=R^Y=R^B$ using Proposition~\ref{prop:r-equality}. We first compare ancestor correlators by matching the graph expansions, then use the $S$-operator to obtain ordinary descendants. A dilaton contraction gives the free energies.

\subsection{The leaf substitution}
The recursion graph sum is naturally written in second-kind differentials $\theta_\alpha^k$. We introduce a second basis $\thh_\alpha^k$ in which ancestor inputs can be substituted directly. The two bases are related by the common matrix $R$.

Define second-kind differentials
\begin{equation}\label{eq:theta}
 \theta_\alpha^k(p)=
 -\frac{(2k-1)!!}{2^k\sqrt{-2}}\,
 \Res_{q=p_\alpha}\frac{B^-(p,q)}{s_\alpha(q)^{2k+1}},
 \qquad k\geq0,\qquad \sqrt{-2}=i\sqrt2.
\end{equation}
Their principal parts at $p_\alpha$ are $-(2k+1)!!\,ds_\alpha/(2^k\sqrt{-2}\,s_\alpha^{2k+2})$. Write $\theta_\alpha(z)=\sum_{k\geq0}\theta_\alpha^kz^k$ and define the forms $\thh_\alpha^k$ by the triangular identity
\begin{equation}\label{eq:leaf-dressing}
 \theta_\alpha(z)=\sum_\beta R_{\beta\alpha}(z)\thh_\beta(z),
 \qquad \thh_\beta(z)=\sum_{k\geq0}\thh_\beta^kz^k.
\end{equation}
For a polynomial $U(z)=\sum_{\alpha,k}U_k^\alpha z^k\eh_\alpha$, let $\mathcal L_U$ be the substitution
\begin{equation}\label{eq:ancestor-dictionary}
 \mathcal L_U(\thh_\alpha^k)=(-1)^{k+1}U_k^\alpha.
\end{equation}
The principal parts make these forms linearly independent. For a fixed stable $(g,n)$ the recursion has a finite expansion in their tensor products, so the substitution is well-defined.

For the shifted CohFT \eqref{eq:primary-shift}, define the $T$-equivariant ancestor correlator by
\[
 \mathcal A_{g,n}^{Y,\tau}(U_1,\ldots,U_n)
 =\int_{\M_{g,n}}
   \Omega^{Y,\tau}_{g,n}\bigl(U_1(\psi_1),\ldots,U_n(\psi_n)\bigr).
\]
Here $\psi_i$ is the cotangent-line class at the $i$th marking on $\M_{g,n}$. The normalized idempotents in this formula are taken at $\tau$. The divisor equation incorporates the divisor part of $\tau$ into the mirror parameters.

\begin{proposition}\label{prop:ancestor}
For $n>0$ and $2g-2+n>0$,
\[
 (\mathcal L_{U_1}\otimes\cdots\otimes\mathcal L_{U_n})\omega_{g,n}
 =\mathcal A_{g,n}^{Y,\tau}(U_1,\ldots,U_n).
\]
\end{proposition}
\begin{proof}
The local graph expansion of topological recursion \cite[Theorems~3.7 and~4.1]{DOSS}, in the reduced normalization of \cite[Section~4.5]{FLM}, has an overall factor $(-1)^{g-1}$, vertex factors $v_\alpha^{2-2g_v-\operatorname{val}(v)}$, ordinary leaves $\theta_\alpha^k$, and edge series
\[
 \frac{I-R(z)^TR(w)}{z+w}.
\]
To identify its translation, expand the odd part of $\yt$:
\begin{equation}\label{eq:odd-y}
 \yt_\alpha(s)-\yt_\alpha(-s)
 =2\sum_{k\geq1}\frac{2^{k-1}}{(2k-1)!!}h_{\alpha,k}s^{2k-1}.
\end{equation}
The metric normalization gives $v_\alpha=h_{\alpha,1}/\sqrt{-2}$. The kernel-period identity for the unit in Proposition~\ref{prop:b-flat}, followed by Gaussian moments, gives
\begin{equation}\label{eq:unit-row}
 \sum_\beta v_\beta R_{\beta\alpha}(z)
 =\frac1{\sqrt{-2}}\sum_{k\geq1}h_{\alpha,k}z^{k-1}.
\end{equation}
Thus the stable dilaton leaf is $-[z^k]\,z\sum_\beta v_\beta R_{\beta\alpha}(z)$ for $k\geq2$. This is the flat-unit compatibility of \rev{\cite[Proposition~6.3]{FLM}}.

Givental--Teleman reconstruction \cite{GiventalQuant,Teleman} uses $R^{-1}(z)=R(-z)^T$. Its edge and translation series are
\[
 \frac{I-R(-z)^TR(-w)}{z+w},
 \qquad zv-zR(-z)^Tv.
\]
For an ordinary leaf, \eqref{eq:leaf-dressing} and \eqref{eq:ancestor-dictionary} give
\[
 \begin{aligned}
 \mathcal L_U(\theta_\alpha^k)
 &=\sum_{\beta}\sum_{j=0}^k
       (R_j)_{\beta\alpha}(-1)^{k-j+1}U_{k-j}^{\beta}\\
 &=(-1)^{k+1}[z^k]\bigl(R(-z)^TU(z)\bigr)_\alpha\\
 &=(-1)^{k+1}[z^k]\bigl(R^{-1}(z)U(z)\bigr)_\alpha.
 \end{aligned}
\]
Here $R_0=I$, and $U(z)$ in the last two lines denotes its column in the normalized frame. The last line is the A-model marking factor with sign $(-1)^{k+1}$. Comparing the edge series gives $(-1)^{k+j+1}$ at half-edge heights $k,j$; comparing the translation series gives $(-1)^{k+1}$ for a dilaton leaf of height $k$.
For a connected graph with $e$ edges, $n$ ordinary leaves and $m$ dilaton leaves, the vertex dimension equations give
\[
 \sum_{\text{half-edges and leaves}} k
       =3g-3-e+n+m.
\]
The product of the signs is therefore $(-1)^{g-1}$, canceling the same overall factor in the local recursion expansion. The Frobenius data and calibration agree by Proposition~\ref{prop:r-equality}, so every graph contribution agrees. Since dilaton heights are at least two, the dimension equation also gives $m+e\leq3g-3+n$. Only finitely many graphs and matrix coefficients contribute.
\end{proof}

\subsection{Ordinary descendants}
The ancestor classes above use cotangent lines after stabilization to $\M_{g,n}$. Ordinary descendants use cotangent lines on the stable-map space. The $S$-operator accounts for the rational tails contracted by stabilization.

For $\gamma_i\in H$ and integers $k_i\geq0$, define the $T$-equivariant descendant Gromov--Witten bracket, for types admitting stable maps, by
\[
 \left\langle\gamma_1\psi_1^{k_1},\ldots,\gamma_n\psi_n^{k_n}
                       \right\rangle_{g,n,\beta}^{Y}
 :=
 \sum_{F\subset\M_{g,n}(Y,\beta)^T}
 \int_{[F]^{\rm vir}}
 \frac{\prod_{i=1}^n\ev_{i,F}^*\gamma_i\,(\psi_i|_F)^{k_i}}
      {e_T(N_F^{\rm vir})}.
\]
The sum runs over connected fixed components, with the virtual classes and normal bundles used in \eqref{eq:cohft}. Here $\psi_i=c_1(\mathbb L_i)$ on $\M_{g,n}(Y,\beta)$, where $\mathbb L_i$ has fiber $T^*_{x_i}C$ at a stable map $[f:C\to Y;x_1,\ldots,x_n]$. Extend the bracket multilinearly to polynomial descendant insertions. The insertion $u/(z-\psi)$ means $\sum_{k\geq0}u\psi^kz^{-k-1}$. Define the $T$-equivariant Gromov--Witten operator $S_\tau(z)$ by
\begin{equation}\label{eq:s-operator}
 \eta(v,S_\tau(z)u)=\eta(v,u)+
 \sum_{\beta,m}^{\rm stable}\frac{q^\beta}{m!}
 \left\langle v,\frac{u}{z-\psi},\tau,\ldots,\tau
                                      \right\rangle_{0,m+2,\beta}^{Y}.
\end{equation}
The sum omits $(\beta,m)=(0,0)$ and includes positive-degree two-point invariants. For a polynomial input $u(z)$, put
\[
 U^\alpha(z)=
   \bigl[\eta(\eh_\alpha,S_\tau(z)u(z))\bigr]_+,\qquad
 U(z)=\sum_\alpha U^\alpha(z)\eh_\alpha,
\]
where $[\ ]_+$ retains nonnegative powers of $z$.

\begin{theorem}\label{thm:descendants}
Under the mirror map \eqref{eq:mirror-map}, with compatible choices of critical points and normalized idempotents, for $n>0$ and $2g-2+n>0$, the polynomial inputs $u_i(z)$ and corresponding $U_i(z)$ defined above satisfy
\begin{equation}\label{eq:descendant-dictionary}
 (\mathcal L_{U_1}\otimes\cdots\otimes\mathcal L_{U_n})\omega_{g,n}
 =
 \sum_{\beta,m}\frac{q^\beta}{m!}
 \left\langle u_1(\psi_1),\ldots,u_n(\psi_n),
                      \tau,\ldots,\tau\right\rangle_{g,n+m,\beta}^{Y}.
\end{equation}
\end{theorem}
\begin{proof}
The ancestor--descendant relation replaces each descendant input $u$ by $[S_\tau u]_+$ \cite[Theorem~5.1 and Proposition~5.3]{GiventalQuant}; it follows by expanding the boundary divisors where a rational tail is contracted when the map and the extra markings are forgotten and the source curve is stabilized \cite[Theorem~1.1]{KM}. Proposition~\ref{prop:ancestor} then gives \eqref{eq:descendant-dictionary}.
\end{proof}
General formal descendant inputs are interpreted coefficientwise in their variables.

The continuation assertion also applies to the descendant leaves. At $\tau=0$ write $S=\sum_{r\geq0}S_r(q)z^{-r}$ with $S_0=I$. In the following operator formulas, $D_i$ also denotes classical cup multiplication by this divisor, and $C_i(q)$ denotes its quantum multiplication operator. For $r\geq1$ one has $S_r(0)=0$, and the divisor equation gives
\[
 q_i\partial_{q_i}S_{r+1}=C_i(q)S_r-S_rD_i,\qquad r\geq0.
\]
Here $D_i=C_i(0)$. The Bryan--Gholampour finite root formula makes $C_i(q)$ analytic near large radius. Writing $C=\sum_iC_i$ and $D=\sum_iD_i$, the recursion
\[
 S_{r+1}(q)=\int_0^1
       \bigl(C(sq)S_r(sq)-S_r(sq)D\bigr)\,\frac{ds}{s}
\]
proves coefficientwise convergence; its numerator has zero constant term. The flat connection then continues each coefficient. For fixed genus and descendant degrees, reconstruction uses finitely many coefficients of $R$ and of these leaves. Its B-model expression therefore gives coefficientwise continuation through adjacent semisimple domains.

\subsection{The free energies}
For $g\geq2$, the free energy can be recovered from the one-point ancestor correlator with insertion $\psi\,1$. We show that the residue defining $F_g^{\TR}$ implements this insertion.

For $g\geq2$ define the $T$-equivariant Gromov--Witten free energy
\[
 F_g^{\GW}(Y;\tau)=
 \sum_{\beta,m}\frac{q^\beta}{m!}
          \langle\tau,\ldots,\tau\rangle_{g,m,\beta}^{Y}
 =\int_{\M_g}\Omega^{Y,\tau}_{g,0}.
\]
Choose local primitives $d\Phi_\alpha=\yt_\alpha\,d\xh$ and use the residue definition of the recursion free energy \cite{EO}:
\begin{equation}\label{eq:tr-free-energy}
 F_g^{\TR}=\frac{1}{2-2g}\sum_\alpha
 \Res_{p=p_\alpha}\Phi_\alpha(p)\omega_{g,1}(p).
\end{equation}
Here $\Phi_\alpha$ takes values in $\mathcal L_{\rm sgn}$, so $\Phi_\alpha\omega_{g,1}$ is a scalar one-form under the canonical pairing $\mathcal L_{\rm sgn}^{\otimes2}\simeq\OO_{Q^\circ}$. Stable recursion forms have no residues, so the choice of additive constants in these primitives does not affect $F_g^{\TR}$.

\begin{theorem}\label{thm:free-energies}
For $g\geq2$, under the mirror map \eqref{eq:mirror-map},
\[
 F_g^{\GW}(Y;\tau)=F_g^{\TR}.
\]
\end{theorem}
\begin{proof}
Define the contraction $\ell(\xi)=-\sum_\alpha\Res_{p_\alpha}\Phi_\alpha\xi$. The principal parts in \eqref{eq:theta}, together with \eqref{eq:odd-y}, give
\[
 \ell(\theta_\alpha^0)=0,\qquad
 \ell(\theta_\alpha^k)=\frac{h_{\alpha,k}}{\sqrt{-2}}
 \quad(k\geq1).
\]
Equations \eqref{eq:leaf-dressing} and \eqref{eq:unit-row} imply
\[
 \ell(\thh_\alpha^k)=\delta_{k1}v_\alpha.
\]
This is exactly the substitution for the ancestor input $z\,1$: the sign in \eqref{eq:ancestor-dictionary} at height one is positive. By Proposition~\ref{prop:ancestor} and the flat-unit axiom of the shifted CohFT,
\[
 \ell(\omega_{g,1})
   =\int_{\M_{g,1}}\psi_1\Omega^{Y,\tau}_{g,1}(1)
   =(2g-2)\int_{\M_g}\Omega^{Y,\tau}_{g,0}.
\]
Since $F_g^{\TR}=\ell(\omega_{g,1})/(2g-2)$, this proves $F_g^{\GW}(Y;\tau)=F_g^{\TR}$.
\end{proof}

\section{Integral cycles for the resolution}\label{sec:cycles}
The descendant dictionary of Section~\ref{sec:correlators} is expressed in cohomological inputs. We now assign global mirror contours to geometric K-classes; Section~\ref{sec:gamma-periods} will identify their oscillatory integrals with Gamma-framed Gromov--Witten series.
Put $r=l-2$ and let $\mathsf u$ denote the character of $T$ of weight $\nu$. Let $K_Z^T(Y)$ be the Grothendieck group of equivariant coherent sheaves on $Y$ supported on the zero section $i:Z\hookrightarrow Y$, equivalently the $K_0$ of their bounded derived category. D\'evissage gives an isomorphism of K-groups
\[
 i_*:K_T(Z)=K_0(D^b\operatorname{Coh}^T(Z))
       \xrightarrow{\sim}K_Z^T(Y).
\]
Indeed, if $\mathcal I_Z$ is the ideal of $Z$, every supported sheaf $\mathcal F$ has a finite filtration $\mathcal F\supset\mathcal I_Z\mathcal F\supset\cdots \supset\mathcal I_Z^N\mathcal F=0$. Its successive quotients are sheaves on $Z$, and their classes sum to $[\mathcal F]$; this construction gives the inverse on K-groups.
In this product geometry the projection $p:Y\to Z$ realizes the inverse by $p_*$: each supported sheaf lies on a finite thickening of $Z$, so its pushforward to $Z$ is coherent and exact, and $p_*i_*=\id$.

There is also a natural map $j:K_Z^T(Y)\to K_T(Y)$ forgetting the support condition. Writing $p:Y\to Z$, the Koszul sequence gives
\[
 j([i_*F])=(1-\mathsf u^2)p^*[F].
\]
The map $j$ is injective in our equivariant setting (Lemma~\ref{lem:k-basis}). Under $\mathsf u\mapsto\mathsf q^{-1}$, the normal-line factor becomes $1-\mathsf q^{-2}$, which reappears in the zero-end capping coefficients below.

\subsection{A geometric K-basis}
Write the binary dihedral group of order $4r$ as
\[
 G=\langle a,b\mid a^{2r}=1,\ b^2=a^r,\ bab^{-1}=a^{-1}\rangle.
\]
Its irreducible representations, in the order used below, are $\rho_0,\rho_1,\rho_2,\ldots,\rho_r,\rho_L,\rho_R$. The four representations $\rho_0,\rho_1,\rho_L,\rho_R$ have dimension one. The remaining $r-1$ representations $\rho_2,\ldots,\rho_r$ have dimension two. Here $\rho_0$ is trivial, $\rho_1(a)=1$, $\rho_1(b)=-1$, and $\rho_v$ has $a$-eigenvalues $e^{\pm\pi i(v-1)/r}$ for $2\le v\le r$. The characters $\rho_L,\rho_R$ satisfy $\rho_L(a)=\rho_R(a)=-1$ and $\rho_L(b)=-\rho_R(b)=1$ for even $r$, or $i$ for odd $r$.

The geometric McKay correspondence was constructed by Gonzalez-Sprinberg--Verdier \cite{GSV}; Ito--Nakamura \cite[Theorem~9.3]{IN} identify $Z$ with $G\text{-Hilb}(\C^2)$. This scheme parametrizes $G$-clusters: $G$-invariant zero-dimensional subschemes $\xi\subset\C^2$ whose coordinate algebra $H^0(\xi,\OO_\xi)$ is isomorphic to the regular representation $\C[G]$ as a $G$-module \cite[Section~1]{BKR}. Such a subscheme has length $4r$. Let $\mathcal R_\rho$ be the tautological bundle with fiber $\operatorname{Hom}_G(\rho,H^0(\xi,\OO_\xi))$ at $[\xi]\in Z$. Its rank is $\dim\rho$, so the four one-dimensional representations give line bundles and the others give rank-two bundles. Coordinate functions have weight $-\nu$. Let $\mathsf q$ be an invertible formal variable for the cycle coefficients. For oscillatory integrals we specialize $\mathsf q=e^{2\pi i\nu/z}$. For equivariant bundles $E,F$ on $Z$, write
\[
 \chi_T(E,F)=\sum_{k\geq0}(-1)^k
       \operatorname{char}_T\operatorname{Ext}_Z^k(E,F).
\]
Here the Ext groups carry their induced $T$-action; their graded characters are interpreted as rational functions of $\mathsf u$.

\begin{lemma}\label{lem:k-basis}
The classes $P_\rho=i_*\mathcal R_\rho^\vee$ form a basis of $K_Z^T(Y)$ over $\mathbb Z[\mathsf u^{\pm1}]$. If $A$ is the adjacency matrix of the affine $D_l$ diagram, in the order above, then the matrix of Euler pairings, with row $\rho$ and column $\sigma$, is
\begin{equation}\label{eq:euler-matrix}
 \begin{aligned}
 \chi(\mathsf q)
 &:=\left(
   \left.\chi_T(\mathcal R_\rho,\mathcal R_\sigma)
                  \right|_{\mathsf u=\mathsf q}\right)_{\rho,\sigma}\\
 &=\big((1+\mathsf q^{-2})I-\mathsf q^{-1}A\big)^{-1}.
 \end{aligned}
\end{equation}
\end{lemma}
\begin{proof}
Let $\mathcal U\subset Z\times\C^2$ be the universal cluster. Write $\pi_Z:Z\times\C^2\to Z$ and $\pi_{\C^2}:Z\times\C^2\to\C^2$ for the projections. The surface derived McKay equivalence of Kapranov--Vasserot \cite[Section~1.4]{KV}, in the universal-cluster formulation of Bridgeland--King--Reid \cite[Theorem~1.2]{BKR}, has inverse kernel $\OO_{\mathcal U}^{\vee}\otimes\pi_{\C^2}^*\omega_{\C^2}[2]$. Here $\OO_{\mathcal U}^{\vee}$ denotes the derived dual on $Z\times\C^2$. Relative duality over $Z$ identifies its pushforward with $((\pi_Z)_*\OO_{\mathcal U})^\vee$. Since $(\pi_Z)_*\OO_{\mathcal U}=\bigoplus_\rho\rho\otimes\mathcal R_\rho$, the inverse sends $\OO_{\C^2}\otimes\rho$ to $\mathcal R_\rho^\vee$.

These constructions are equivariant for the commuting scalar action. Equivariant homotopy invariance \cite{Thomason} on $\C^2$ therefore gives the basis $[\mathcal R_\rho^\vee]$ of $K_T(Z)$ over $\mathbb Z[\mathsf u^{\pm1}]$. Applying d\'evissage gives the asserted basis $P_\rho$ of $K_Z^T(Y)$. Since $p^*:K_T(Z)\to K_T(Y)$ is an isomorphism and multiplication by $1-\mathsf u^2$ is injective on this free module, the Koszul identity above also shows that the forget-support map $j$ is injective.

For the Euler matrix, put $V=(\C^2)^\vee$, the space of linear coordinate functions, and set
\[
 (H_d)_{\rho\sigma}
 =\dim\operatorname{Hom}_G
       (\rho,\operatorname{Sym}^dV\otimes\sigma).
\]
The derived equivalence computes the graded Euler matrix of the dual bundles as $\sum_{d\ge0}\mathsf u^{-d}H_d$. As $G$-representations, $V\simeq\C^2$ and $\det V$ is trivial. The Clebsch--Gordan relation therefore gives
\[
 H_0=I,\qquad H_1=A,\qquad
 AH_d=H_{d+1}+H_{d-1}\quad(d\ge1).
\]
Consequently
\[
 \sum_{d\ge0}s^dH_d=((1+s^2)I-sA)^{-1}.
\]
Taking $s=\mathsf u^{-1}$, which records the coordinate degree, proves \eqref{eq:euler-matrix}.
This matrix is symmetric, and $\chi_T(\mathcal R_\rho,\mathcal R_\sigma) =\chi_T(\mathcal R_\sigma^\vee,\mathcal R_\rho^\vee)$, which proves the stated convention.
\end{proof}

\subsection{The integral homology group}
Let $\mathscr R=\mathbb Z[\mathsf q^{\pm1}]$. On the quotient $w$-sphere remove the zeros $c_1,\ldots,c_l$ of $\lambda$ and its poles $2,-2,\infty$. Let $M$ be the real oriented blowup at these punctures and put $\partial_PM=\partial_2M\sqcup\partial_{-2}M\sqcup\partial_\infty M$. Define $\mathcal L$ to be the locally constant sheaf of free rank-one $\mathscr R$-modules on $M$ with counterclockwise monodromies
\begin{equation}\label{eq:cycle-monodromy}
 \mathsf q^{-2}\quad\text{at }c_j,\qquad
 -\mathsf q^2\quad\text{at }2,-2,\qquad
 \mathsf q^{2r}\quad\text{at }\infty.
\end{equation}
Their product is one. After $\mathsf q=e^{2\pi i\nu/z}$, on the punctured curve $\mathcal L$ is the tensor product of the sign local system $\mathcal L_{\rm sgn}$ from Section~\ref{subsec:equivariant-recursion} and the local system of branches of $e^{-W/z}=\lambda^{-2\nu/z}$. The former supplies the signs at $\pm2$; the latter supplies the powers of $\mathsf q$. The cycle correspondence uses $\mathsf u\mapsto\mathsf q^{-1}$, as the modified Chern character sends a weight-$\nu$ line to $e^{-2\pi i\nu/z}$. Our cycle group is
\begin{equation}\label{eq:mixed-homology}
 \mathcal H_{\mathscr R}^+=H_1(M,\partial_PM;\mathcal L).
\end{equation}
After $\mathsf q=e^{2\pi i\nu/z}$, the oscillatory integrals have integrand
\begin{equation}\label{eq:period-integrand}
 e^{-W/z}\frac{d\mu}{\mu}
 =\lambda^{-2\nu/z}\frac{dw}{\sqrt{w^2-4}}.
\end{equation}
We integrate along one lift to the $\mu$-sphere. The local-system coefficients record the branches of the multivalued integrand, as in the twisted-homology formulation of hypergeometric integrals \cite[Chapter~2]{AK}. The allowed ends are the poles of $\lambda$; for $\Re(2\nu/z)>0$ they are decay ends.

\rev{We consider the radial period paths $[0,\kappa_j^{-1}]$ and $[0,1]$ that appear in Ma's period basis \cite[Propositions~3.4.12--3.4.13]{MaThesis}.} We now specify the integral coefficients, branches, and cappings of the contours used here.
Choose first $-2<c_1<\cdots<c_l<2$ and arcs with disjoint interiors in the lower half-plane from infinity to the punctures. For $0\leq j\leq l-1$, denote the arc to $c_{j+1}$ by $\gamma_j$; denote the arcs to $2,-2$ by $e_+,e_-$. These are the images of the radial segments from $0$ to $\kappa_{j+1}^{-1}$ in the upper $\mu$-half-disc, together with the segments from $0$ to $1$ and $-1$. Each arc carries a flat section of $\mathcal L$, specified by its lift and a branch of $\lambda^{-2\nu/z}$. Use the branch $\lambda^{-2\nu/z}\sim K^{-2\nu/z}\mu^{2r\nu/z}$ near $\mu=0$, continued with $0\le\arg\mu\le\pi$. This fixes the sections on $\gamma_j$ and $e_+$. On $e_-$ multiply the continued section by $\mathsf q^{-2r}$; its leading coefficient along $\mu=-s$, $s>0$, is $K^{-2\nu/z}\mathsf q^{-r}$. Continue these paths and their sections together.
The paths $\gamma_j$ have zero ends. To obtain classes in $\mathcal H_{\mathscr R}^+$, truncate such an end and attach a small counterclockwise loop $\ell_j$, followed by the return path with its continued section. Define the capped cycle by
\[
 \delta_j
 =\mathsf q^2(\gamma_j+\ell_j-\mathsf q^{-2}\gamma_j),
 \qquad 0\le j\le l-1.
\]
The two radial parts have opposite orientations, and the loop monodromy $\mathsf q^{-2}$ cancels their boundary at the zero. In particular, $\delta_0$ is the capped contour around $c_1$. Under the natural map
\[
 \iota_*:\mathcal H_{\mathscr R}^+
       \longrightarrow H_1(M,\partial M;\mathcal L),
\]
the loop is killed and
\[
 \iota_*(\delta_j)=(\mathsf q^2-1)[\gamma_j].
\]
Thus the scalar multiple of the zero arc describes the image in the full-boundary relative group; the capped chain defines the integral class in $\mathcal H_{\mathscr R}^+$.
Set $f=\mathsf q^2-1$ and $v=\mathsf q^2+1$.

\begin{lemma}\label{lem:cycle-basis}
The group \eqref{eq:mixed-homology} is free of rank $l+1$, with basis
\[
 \delta_1,\ldots,\delta_{l-1},e_+,e_-.
\]
The generators $\delta_0,\ldots,\delta_{l-1},e_+,e_-$ satisfy the single relation
\begin{equation}\label{eq:arc-relation}
 \sum_{j=0}^{l-1}\mathsf q^{2j}\delta_j
 -\mathsf q^{2l-2}v e_++\mathsf q^{2l-4}v e_-=0.
\end{equation}
\end{lemma}
\begin{proof}
Cut along these arcs and read the boundary clockwise. For puncture monodromies $t_j$, its coefficient is $(\prod_{k<j}t_k^{-1})(t_j^{-1}-1)$. This gives \eqref{eq:arc-relation}, with $f\gamma_j$ in place of $\delta_j$, for homology relative to all boundary circles. To impose the zero-end condition, put $Z_0=\bigcup_j\partial_{c_j}M$. Poincar\'e--Lefschetz duality and the relative cohomology sequence identify \eqref{eq:mixed-homology} with the kernel of
\[
 H^1(M;\mathcal L)\longrightarrow H^1(Z_0;\mathcal L)
       =\bigoplus_j\mathscr R/(f).
\]
Indeed $H^0(Z_0;\mathcal L)=0$, since $f$ is a nonzero divisor. Restriction reads each zero-arc coefficient modulo $f$, up to a unit. Its kernel is therefore generated by the contours $\delta_j$ around $c_{j+1}$, for $0\leq j\leq l-1$, and the arcs $e_\pm$ from $\infty$ to $\pm2$, subject to the displayed relation. Its coefficient of $\delta_0$ is one, so the relation expresses $\delta_0$ in terms of the remaining generators, which form the asserted free basis.
\end{proof}

\begin{theorem}\label{thm:k-cycles}
There is a semilinear isomorphism
\[
 \operatorname{mir}_Z:K_Z^T(Y)\longrightarrow\mathcal H_{\mathscr R}^+,
 \qquad \operatorname{mir}_Z(\mathsf u F)=\mathsf q^{-1}\operatorname{mir}_Z(F).
\]
Writing $M_\rho=\operatorname{mir}_Z(P_\rho)$, it is given by
\begin{align}\label{eq:explicit-cycles}
 M_0&=\mathsf q^{r-1}e_-,\notag\\
 M_1&=-\sum_{j=1}^{r+1}\mathsf q^{2j-r-3}\delta_j
       +\mathsf q^{r-1}(1+\mathsf q^2)e_+-\mathsf q^{r-1}e_-,\notag\\
 M_a&=-\sum_{j=a}^{r+1}\mathsf q^{2j-a-r-2}\delta_j
       +\mathsf q^{r-a}(1+\mathsf q^2)e_+\quad(2\le a\le r),\\
 M_L&=-\mathsf q^{-1}\delta_{r+1}+\mathsf q e_+,
 \qquad M_R=\mathsf q^{-1}e_+.\notag
\end{align}
The cycles continue by the Gauss--Manin connection on the smooth mirror family.
\end{theorem}
\begin{proof}
All coefficients are integral Laurent polynomials. Conversely,
\begin{align*}
 e_+&=\mathsf q M_R,&e_-&=\mathsf q^{1-r}M_0,\\
 \delta_{r+1}&=\mathsf q^3M_R-\mathsf q M_L,&
 \delta_r&=\mathsf q^2\bigl(\mathsf q(M_L+M_R)-M_r\bigr),\\
 \delta_a&=\mathsf q^{r+2-a}(\mathsf q M_{a+1}-M_a)\quad(2\le a<r),&
 \delta_1&=\mathsf q^{r+1}(\mathsf q M_2-M_1-M_0).
\end{align*}
Thus the images form a basis over $\mathscr R$. Lemma~\ref{lem:k-basis} proves the isomorphism, and transporting the contours with their flat sections gives its continuation. Section~\ref{sec:gamma-periods} identifies the corresponding oscillatory integrals with the resolution's equivariant Gromov--Witten central charges.
\end{proof}

\subsection{Line bundles and their contours}
Let $\mathcal L_i=(\det\mathcal R_{\rho_i})^\vee$ for $i=1,\ldots,r,L,R$, with $E_L=E_{l-1}$ and $E_R=E_l$. The geometric McKay correspondence \cite[Theorem~2.2]{GSV} gives degree $-1$ on $E_i$ and degree zero on every other exceptional curve. We use the equivariant structures inherited from the universal cluster. Their fixed-component restrictions, computed below, verify $c_1^T(\mathcal L_i)=-D_i$, with $D_i=-\sum_j(C_{D_l}^{-1})_{ji}[E_j]$ the class dual to $E_i$. \rev{The restrictions and tautological-to-line-bundle decomposition agree with the formulas of Brini--Ma \cite[Proposition~5.4 and Lemma~5.6]{BM}, after translating the determinant conventions. The cluster calculation below relates these classes to the integral contours.}
\begin{lemma}\label{lem:bundle-restrictions}
The restrictions of the fundamental line bundles are
\begin{equation}\label{eq:bundle-restrictions}
\begin{array}{c|c|c|c|c}
 &p_k,\ 1\le k<r&E&p_L&p_R\\ \hline
 \mathcal L_a,\ 1\le a\le r
 &\begin{cases}\mathsf u^{2a},&a<k,\\
                 \mathsf u^{2r},&a\ge k\end{cases}
 &\mathsf u^{2a}\OO_E(-\delta_{ar})&\mathsf u^{2a}&\mathsf u^{2a}\\
 \mathcal L_L&\mathsf u^r&\mathsf u^r\OO_E&\mathsf u^{r+2}&\mathsf u^r\\
 \mathcal L_R&\mathsf u^r&\mathsf u^r\OO_E&\mathsf u^r&\mathsf u^{r+2}
\end{array}
\end{equation}
For $2\le a\le r$, the tautological classes satisfy
\begin{equation}\label{eq:tautological-line-decomposition}
 [\mathcal R_{\rho_a}]
 =\mathsf u^{a-1}[\mathcal L_a^{-1}]
                       +\mathsf u^{1-a}[\OO_Z].
\end{equation}
\end{lemma}
\begin{proof}
Use the type-$D$ cluster ideals of \cite[Section~13]{IN}, with coordinate functions of weight $-\nu$. Choose $a(x,y)=(e^{\pi i/r}x,e^{-\pi i/r}y)$ and $b(x,y)=(y,-x)$. The ideal at $p_1$ is $(xy,x^{2r}+y^{2r})$, and
\[
 (x^ky,xy^k,x^2y^2,x^{2r+2-k},y^{2r+2-k}),
 \qquad 2\le k\le r,
\]
represents $p_k$ for $k<r$ and the point $[0:1]$ of $E$ for $k=r$. These homogeneous ideals are preserved by scalar multiplication. The fixed curve is the projective family
\begin{equation}\label{eq:fixed-cluster-family}
 \bigl(x^2y^2,\ ux^{r+1}-vxy^r,\
 uy^{r+1}-(-1)^rvx^ry,\
 x^{r+1}y,\ xy^{r+1},\ x^{r+2},\ y^{r+2}\bigr),
 \qquad [u:v]\in\PP^1.
\end{equation}
The two varying generators are exchanged by $b$, up to its prescribed sign. Modulo $x^2y^2$, the dimensions below degree $r+1$ are $1,2,3,4,\ldots,4$; the two independent relations leave dimension two in degree $r+1$, and higher degrees vanish.
These constant graded dimensions make the family finite flat of length $1+2+3+4(r-2)+2=4r$. To verify its $G$-representation, specialize to $[u:v]=[1:0]$. The ideal becomes $(x^2y^2,x^{r+1},y^{r+1})$, and a monomial basis of the quotient is
\[
 x^iy^j,\qquad 0\le i,j\le r,\quad \min(i,j)\le1.
\]
For $a^m\ne1$, put $\zeta=e^{\pi im/r}$. Summing the diagonal coefficients gives
\[
 \operatorname{Tr}(a^m)
 =\frac{(1+\zeta)(\zeta^{-r}-\zeta^r)}{1-\zeta}=0.
\]
For $ba^m$ only $1$ and $xy$ contribute to the trace, with values $1$ and $-1$. The character is therefore $4r$ at the identity and zero elsewhere, the regular character. In a finite flat $G$-family the isotypic summands have locally constant rank, so every fiber has this representation.
Its relation subspace recovers $[u:v]$, identifying the family with $E$. The terminal ideals are $(x^2y^2,x^r\pm\zeta_r^{-1}y^r)$, where $\zeta_r=1$ for even $r$ and $i$ for odd $r$; their quotients also have length $4r$.

The first copy of $\rho_a$, $2\le a\le r$, occurs in degree $a-1$. At $p_k$ its second copy occurs in degree $2r+1-a$ if $k\le a$, and in degree $a+1$ if $k>a$. On $E$ and at the terminal points the second degree is $a+1$. For $a=r$ the varying multiplicity space on $E$ is the quotient of a trivial rank-two bundle by its tautological line, hence is $\OO_E(1)$; all other multiplicity lines are constant. The character $\rho_1$ occurs in degree $2r$ at $p_1$ and degree two elsewhere. Each terminal character occurs in degree $r$, except at its corresponding terminal point, where its degree is $r+2$. Taking dual determinants gives \eqref{eq:bundle-restrictions}.

For the equivariant Chern classes, put $H=c_1(\OO_E(1))$ and $C=C_{D_l}$. The divisor restrictions give
\[
 -D_i|_E=2\nu(C^{-1})_{ri}-\delta_{ri}H,
\]
where $(C^{-1})_{ra}=a$ and $(C^{-1})_{r,L}=(C^{-1})_{r,R}=r/2$. At the isolated fixed points they give
\[
 \begin{aligned}
 -D_i|_{p_k}
 &=2\nu\bigl((r-k+1)(C^{-1})_{ki}
                  -(r-k)(C^{-1})_{k-1,i}\bigr),\\
 -D_i|_{p_S}&=4\nu(C^{-1})_{Si}\qquad(S=L,R),
 \end{aligned}
\]
with row zero interpreted as zero. These expressions agree with the characters and degree in \eqref{eq:bundle-restrictions}, verifying the chosen equivariant lift.
The same multiplicity calculation proves \eqref{eq:tautological-line-decomposition} after restriction to every fixed component. Equivariant localization \cite{ThomasonLefschetz} has torsion kernel, whereas Lemma~\ref{lem:k-basis} makes $K_T(Z)$ free over $\mathbb Z[\mathsf u^{\pm1}]$, so restriction is injective and proves the identity in $K_T(Z)$.
\end{proof}

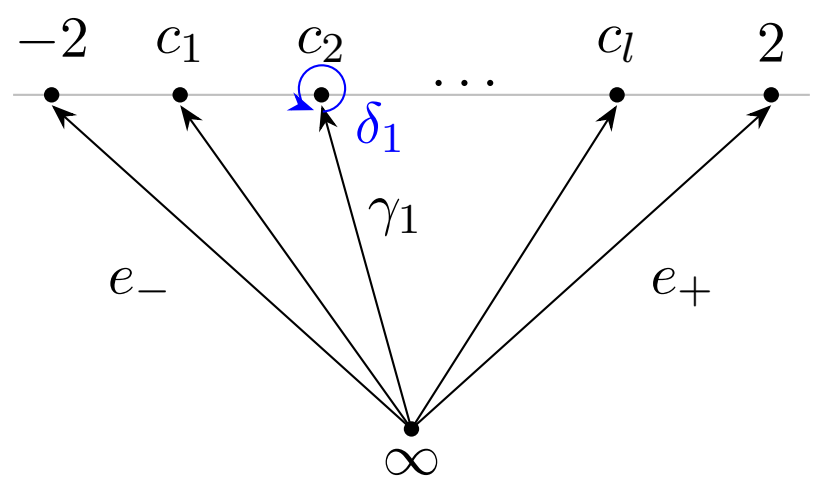
\begin{figure}[ht]
\centering
\begin{tikzpicture}[scale=.85,>=Stealth]
 \coordinate (o) at (0,-2.6);
 \draw[gray!50] (-3.1,0)--(3.1,0);
 \foreach \x/\lab in {-2.8/{-2},-1.8/{c_1},-.7/{c_2},1.6/{c_l},2.8/{2}} {
  \draw[->] (o)--(\x,-.08);
  \fill (\x,0) circle (1.7pt);
  \node[above] at (\x,0.1) {$\lab$};
 }
 \node at (.45,.08) {$\cdots$};
 \fill (o) circle (1.7pt); \node[below] at (o) {$\infty$};
 \node[left] at (-1.7,-1.5) {$e_-$};
 \node[right] at (1.7,-1.5) {$e_+$};
 \node[right] at (-.5,-.95) {$\gamma_1$};
 \draw[blue,->] (-.69,-.13) arc[start angle=-88,end angle=250,radius=.18];
 \node[blue,right] at (-.6,-.25) {$\delta_1$};
\end{tikzpicture}
\caption{On the quotient Toda curve, the arcs from infinity determine the global line-bundle cycles in \eqref{eq:line-bundle-cycles}; a capped contour replaces a zero end by a loop and the return path with its continued section.}
\label{fig:cycles}
\end{figure}

\begin{corollary}\label{cor:line-bundle-cycles}
The supported classes of $\OO_Z,\mathcal L_1,\ldots,\mathcal L_r, \mathcal L_L,\mathcal L_R$ form an integral basis. Their geometric cycles are
\begin{align}\label{eq:line-bundle-cycles}
 i_*\OO_Z&\longmapsto\mathsf q^{r-1}e_-,\notag\\
 i_*\mathcal L_a&\longmapsto
 -\sum_{j=a}^{r+1}\mathsf q^{2j-2a-r-1}\delta_j
 +\mathsf q^{r+1-2a}\bigl((1+\mathsf q^2)e_+-e_-\bigr)
       \quad(1\le a\le r),\\
 i_*\mathcal L_L&\longmapsto-\mathsf q^{-1}\delta_{r+1}
                                     +\mathsf q e_+,\notag\\
 i_*\mathcal L_R&\longmapsto\mathsf q^{-1}e_+.
 \notag
\end{align}
A character twist $\mathsf u^k$ multiplies the contour by $\mathsf q^{-k}$. Here $e_\pm$ run from the pole $\infty$ to the pole $\pm2$ on the quotient sphere, and $\delta_j$ goes around the zero $c_{j+1}$ and returns to $\infty$, with the flat section specified above. See Figure~\ref{fig:cycles}.
\end{corollary}
\begin{proof}
Lemma~\ref{lem:bundle-restrictions} gives \eqref{eq:tautological-line-decomposition}. The one-dimensional representations already give line bundles. Dualizing yields $P_a=\mathsf u^{1-a}i_*\mathcal L_a+\mathsf u^{a-1}P_0$. This is a triangular integral change of basis with invertible diagonal entries. Semilinearity and \eqref{eq:explicit-cycles} now give \eqref{eq:line-bundle-cycles}.
\end{proof}

The positive-degree bundles are explicit as well. Write $N_i$ for the cycles in \eqref{eq:line-bundle-cycles} and $N_0=\mathsf q^{r-1}e_-$.
For $1\le a\le r$, the fixed-component restrictions imply the quadratic K-relation
\[
 ([\mathcal L_a]-\mathsf u^{2a}[\OO_Z])
 ([\mathcal L_a]-\mathsf u^{2r}[\OO_Z])=0.
\]
On $E$ this also uses $([\OO_E(1)]-1)^2=0$. Multiplication by $[\mathcal L_a^{-1}]$ gives
\[
 [\mathcal L_a^{-1}]
 =(\mathsf u^{-2a}+\mathsf u^{-2r})[\OO_Z]
       -\mathsf u^{-2a-2r}[\mathcal L_a].
\]
Applying the cycle correspondence with $\mathsf u\mapsto\mathsf q^{-1}$ gives
\begin{equation}\label{eq:positive-line-cycles}
 i_*\mathcal L_a^{-1}\longmapsto
   (\mathsf q^{2a}+\mathsf q^{2r})N_0-\mathsf q^{2a+2r}N_a.
\end{equation}
For either short arm $i=L,R$, the corresponding expression is $(\mathsf q^r+\mathsf q^{r+2})N_0-\mathsf q^{2r+2}N_i$. Its quadratic K-relation has roots $\mathsf u^r$ and $\mathsf u^{r+2}$.

More generally, every equivariant line bundle has a unique expression
\[
 \mathcal L=\mathsf u^b\prod_{a=1}^r\mathcal L_a^{d_a}
                        \mathcal L_L^{d_L}\mathcal L_R^{d_R}.
\]
Its degree on $E_i$ is $-d_i$. The following Laurent-polynomial rule gives its contour using the same geometric cycles. For a degree-zero line bundle, the fiber characters at the two ends of each nonfixed exceptional $\PP^1$ agree. The exceptional tree is connected, so all fixed components have the same character; on the fixed curve $E$ the underlying degree-zero line bundle is trivial. Injectivity of K-theory restriction, followed by the determinant map, identifies it with that character line. This also proves the uniqueness of the expression above. Set
\[
 w_E=b+2\sum_{a=1}^r a d_a+r(d_L+d_R),\qquad
 w_k=b+2\sum_{a<k}a d_a+2r\sum_{a\ge k}d_a+r(d_L+d_R),
 \quad 1\le k\le r,
\]
so $w_r=w_E$, and put $[d]_u=(u^d-1)/(u-1)$ for any integer $d$. Here $[d]_u$ is a Laurent polynomial, also for negative $d$; for example $[-2]_u=-u^{-1}-u^{-2}$. Define
\begin{align}\label{eq:all-line-bundle-coefficients}
 c_a(\mathsf u)&=\mathsf u^{w_{a+1}-2a}[d_a]_{\mathsf u^{2(r-a)}}\quad(a<r),&
 c_r(\mathsf u)&=d_r\mathsf u^{w_E-2r},\notag\\
 c_L(\mathsf u)&=\mathsf u^{w_E-r}[d_L]_{\mathsf u^2},&
 c_R(\mathsf u)&=\mathsf u^{w_E-r}[d_R]_{\mathsf u^2},\\
 c_0(\mathsf u)&=\mathsf u^{w_E}-\sum_{a=1}^r c_a(\mathsf u)\mathsf u^{2a}
                            -(c_L(\mathsf u)+c_R(\mathsf u))\mathsf u^r.\notag
\end{align}
Then
\begin{equation}\label{eq:all-line-bundle-cycles}
 \operatorname{mir}_Z(i_*\mathcal L)
   =c_0(\mathsf q^{-1})N_0+
                         \sum_{i=1,\ldots,r,L,R}c_i(\mathsf q^{-1})N_i.
\end{equation}
To verify the rule, consecutive chain restrictions differ by $\mathsf u^{w_a}-\mathsf u^{w_{a+1}}$; only $\mathcal L_a$ has a difference there, equal to $\mathsf u^{2r}-\mathsf u^{2a}$. This gives $c_a$ for $a<r$. The degree on $E$ gives $c_r$, the terminal-point restrictions give $c_L,c_R$, and the rank on $E$ gives $c_0$. Integral injectivity of restriction proves the equality of K-classes. In particular, the formula uses integral contours for all positive and negative degrees.

\section{Gamma classes and oscillatory integrals}\label{sec:gamma-periods}
Ma develops an approach to the type-$D$ period comparison through the quantum differential equation and large-radius asymptotics in \cite[Sections~3.4--3.5 and Appendix~A]{MaThesis}. \rev{Here we determine the supported integral contours and their Gamma normalization from component integrals at the partial large-radius boundary. The geometric central-charge interpretation is developed in \cite[Proposition~4.4 and Theorem~5.3]{Fang}.}
The comparison has three steps: the periods and central charges satisfy the same quantum differential equation; their point and rank-two curve boundary data agree; these data determine the global solution. The contour projections will connect the local comparisons to the integral cycles of Section~\ref{sec:cycles}.

We use divisor coordinates for the $T$-equivariant Gromov--Witten series: the separate Novikov symbols are specialized to one, and $\tau=t_0\,1+\sum_i(\log q_i)D_i$, with chosen logarithms. Expansion at large radius is defined by the divisor equation. The normalized genus-zero $T$-equivariant Gromov--Witten $J$-function $J_Y^{\rm norm}(\tau,z)$ has coefficients in $H=H_T^*(Y)\otimes_{\C[\nu]}\C(\nu)$ and degree-zero term $e^{\tau/z}$. The $T$-equivariant operator $S_\tau(z)$ defined in \eqref{eq:s-operator} has degree-zero term $a\mapsto e^{\tau/z}\cup a$. Both include the divisor factor $\exp\bigl(\sum_i(\log q_i)D_i/z\bigr)$. Following the Gamma integral structure of \cite[Section~2.4]{Iritani}, for a class $F\in K_T(Y)$ use the smooth-target convention of \cite[Definition~4.1]{FLYZ}:
\begin{equation}\label{eq:gamma-framing}
 \begin{aligned}
 \kappa_z(F)&=z^3\prod_{a=1}^3\Gamma(1+\delta_a/z)
                 \widetilde{\operatorname{ch}}_z(F)
                 \in z^3H[[z^{-1}]],\\
 \widetilde{\operatorname{ch}}_z(F)&=\sum_b e^{-2\pi i\xi_b/z}
                 \in H[[z^{-1}]].
 \end{aligned}
\end{equation}
Here $\delta_a$ are the equivariant Chern roots of $TY$; for a bundle, $\xi_b$ are its equivariant Chern roots, and the second formula extends additively to K-theory. The notation $\kappa_z(F)$ denotes this Gamma framing; the pointed moduli-space classes used earlier are the $\kappa_j$. The Gamma factors are expanded at $1$ and the exponentials in $z^{-1}$. Put $\hbar=2\nu/z$ and $\mathsf q=e^{\pi i\hbar}$. The variable $\mathsf q$ controls contour monodromy, whereas $q_i$ and $Q=q_r$ are the quantum parameters. For a supported class the Koszul sequence gives
\[
 \widetilde{\operatorname{ch}}_z(i_*E)
 =(1-\mathsf q^{-2})\widetilde{\operatorname{ch}}_z(E).
\]
The normalized $J$-function has leading term $1+\tau/z$. Using the localized equivariant pairing $\eta$, define the central charge in the convention of \cite[Definition~4.2]{FLYZ}:
\begin{equation}\label{eq:gamma-central-charge}
 \begin{aligned}
 Z_F(\tau,z)&=\eta\bigl(\kappa_z(F),J_Y^{\rm norm}(\tau,-z)\bigr)\\
 &\in z^3\C(\nu)[t_0,\log q_1,\ldots,\log q_l]
             [[q_1,\ldots,q_l,z^{-1}]].
 \end{aligned}
\end{equation}
The scalar series is expanded at large radius and $z=\infty$, with the logarithms chosen above. Its one-point interpretation uses the insertion $\kappa_z(F)/(z(z+\psi))$ and the usual unstable completion. For $F\in K_Z^T(Y)$, we write $\kappa_z(F)$ and $Z_F$ for $\kappa_z(j(F))$ and $Z_{j(F)}$.

For the cycles $\operatorname{mir}_Z(F)$, $\int^{\rm reg}$ denotes meromorphic continuation in $\hbar$ from $\Re\hbar>0$, where the integral converges, with the mirror curve and logarithmic branches fixed.
\begin{theorem}\label{thm:gamma-period}
In the large-radius chamber with the branches specified below, every $F\in K_Z^T(Y)$ satisfies
\begin{equation}\label{eq:gamma-period}
 Z_F(\tau,z)=-\frac{\pi iz}{\nu}
    \int_{\operatorname{mir}_Z(F)}^{\rm reg}
           e^{-W/z}\,\frac{d\mu}{\mu}.
\end{equation}
The identity continues as the contours and their flat sections are transported together. It is an identity of single-valued meromorphic functions of $\hbar=2\nu/z\in\C$. For fixed $\nu\ne0$, this gives a meromorphic identity on $z\in\C^*$.
\end{theorem}

\subsection{Contours and the quantum differential equation}\label{subsec:scalar-integrals}
\rev{The continuation sector below agrees with that of \cite[Section~4.1]{BM}. We specify the branches for $\lambda^{-\hbar}$ and the sign local system explicitly.}
Continue the arcs of Section~\ref{sec:cycles} by moving $\kappa_j^{-1}$ from the upper unit semicircle into the upper half-disc, preserving their angular order and transporting the flat sections. In a large-radius sector take
\begin{equation}\label{eq:period-lr-sector}
 |\kappa_1|\gg\cdots\gg|\kappa_l|\gg1,
 \qquad \arg\kappa_1<\cdots<\arg\kappa_l<0,
\end{equation}
with the arguments in a small interval. The lifts are the segments $\mu=x/\kappa_j$, $0\le x\le1$, and $e_+$ lifts to $[0,1]$. Near $\mu=0$ use $\lambda^{-\hbar}\sim K^{-\hbar}\mu^{r\hbar}$, with the logarithms fixed by this continuation. On the negative segment $e_-$ its prescribed section has coefficient $\mathsf q^{-r}$, as in Section~\ref{sec:cycles}. Taking the arguments to zero makes the positive paths pass below the earlier positive zeros. These conventions also specify the partial-boundary family \eqref{eq:arc} with $Q\in(0,1)$, by retaining $\kappa_{r+1}/\kappa_r=Q$ and continuing in $Q$.

For a transported cycle $\Gamma$ define
\begin{equation}\label{eq:resolution-periods}
 I_\Gamma=\int_\Gamma^{\rm reg}\lambda^{-\hbar}\frac{d\mu}{\mu}.
\end{equation}
The capped zero contours are compact, and the remaining ends lie at poles, where the integral converges for $\Re\hbar>0$. Meromorphic continuation in $\hbar$ follows by Taylor subtraction at the ends. Indeed, if $\lambda=x^m a(x)$ with $a(0)\ne0$, each subtracted term integrates to $A_k(\hbar)\epsilon^{-m\hbar+b+k}/(-m\hbar+b+k)$, where $\epsilon>0$ is a small cutoff, $b=1+\operatorname{ord}(d\mu/\mu)$, and a branch of $\log a$ is fixed. The remainder is holomorphic on successively larger strips.

\begin{proposition}\label{prop:scalar-period-qde}
Let $t_a$ be flat coordinates of the resolution Frobenius structure, with $t_0=2\nu\log K$, and let $C_{ab}^{c}$ be its multiplication coefficients. Both $I_\Gamma$ and $Z_F$ satisfy
\begin{equation}\label{eq:scalar-period-qde}
 -z\partial_a\partial_b I=\sum_c C_{ab}^{c}\partial_c I,
 \qquad -z\partial_0 I=I.
\end{equation}
\end{proposition}
\begin{proof}
\rev{The twisted-period quantum equation is established in \cite[Corollary~3.4.10 and Proposition~3.4.13]{MaThesis}. Our exponent $-\hbar$ corresponds to replacing $z$ by $-z$ in that convention. We give a direct calculation in the normalization of \eqref{eq:scalar-period-qde}.}
For the integral, use $\tau_i=\log\kappa_i$ and $t_0$, which differ from the divisor flat coordinates by a constant linear transformation. Put $u=\log\mu$, $w=\mu+\mu^{-1}$, $v=\mu-\mu^{-1}$, $c_i=\kappa_i+\kappa_i^{-1}$, and $v_i=\kappa_i-\kappa_i^{-1}$. The primary functions are
\[
 f_i=\partial_{\tau_i}W=-\frac{2\nu v_i}{w-c_i},
 \qquad f_0=1,\qquad
 W_u=2\nu v\left(\sum_j\frac1{w-c_j}-\frac{2w}{w^2-4}\right).
\]
For $i\ne k$, partial fractions express $f_if_k$ in the span of $f_i,f_k$, and $\partial_{\tau_i}f_k=0$. For $i=k$, set $g_i=2\nu v/(w-c_i)$. Direct differentiation gives
\[
 \partial_u g_i
 =2\nu\left(-\frac{c_i}{w-c_i}
              -\frac{c_i^2-4}{(w-c_i)^2}\right)
 =\partial_{\tau_i}f_i.
\]
The double pole of $f_i^2-W_ug_i$ cancels. This function has only simple poles at the $c_j$ and is bounded at infinity. Its values at the critical orbits, together with the genus-zero mirror theorem, give
\[
 f_i^2-W_ug_i=\sum_c C_{ii}^{c}f_c.
\]
It follows that
\[
 -z\partial_{\tau_i}^2I_\Gamma-\sum_cC_{ii}^{c}\partial_cI_\Gamma
 =\int_\Gamma^{\rm reg}d(e^{-W/z}g_i)=0.
\]
For generic $\hbar$, replace all puncture ends by compact twisted contours; an exact differential integrates to zero on these contours. Their integrals agree with the convergent end integrals and then with their meromorphic continuations. The unit equation follows from $\partial_0W=1$. The same equations for $Z_F$ are the scalar quantum equations for the $T$-equivariant $J$-function at argument $-z$.
\end{proof}

\subsection{Specialization to the component curves}\label{subsec:gamma-specialization}
Use the family \eqref{eq:arc}, with positive amplitudes and integral exponents
\[
 d_1>\cdots>d_{r-1}>d_f>d_l>0,\qquad r=l-2.
\]
The parameter $Q=\kappa_{r+1}/\kappa_r$ varies in a simply connected neighborhood of a subinterval of $(0,1)$. \rev{The following scalar factors agree, up to the common factor $K$, with the fixed-component monomials in \cite[Proposition~4.12]{BM}:}
\begin{equation}\label{eq:gamma-component-units}
\begin{aligned}
 B_k&=K(\kappa_1\cdots\kappa_{k-1})\kappa_k^{r+1-k}
                                      &&(1\le k<r),\\
 B_f&=K\kappa_1\cdots\kappa_r,&
 B_L&=K\frac{\kappa_1\cdots\kappa_{r+1}}{\kappa_l},
 &B_R&=K\kappa_1\cdots\kappa_l .
\end{aligned}
\end{equation}
Here $B_L$ belongs to $C_l$ and $B_R$ to $C_0$. Write $B_j=a_j(Q)t^{-b_j}$, with the positive branch on the real family. Their exponents are
\begin{equation}\label{eq:gamma-normal-exponents}
\begin{aligned}
 b_k&=\sum_{i<k}d_i+(r+1-k)d_k,&
 b_f&=\sum_{i<r}d_i+d_f,\\
 b_L&=\sum_{i<r}d_i+2d_f-d_l,&
 b_R&=\sum_{i<r}d_i+2d_f+d_l .
\end{aligned}
\end{equation}
Choose the $d_i$ so that these numbers are pairwise distinct. The signs of the component functions are accounted for by the flat sections of the contours.

\begin{lemma}\label{lem:gamma-component-splitting}
For generic $\mathsf q$, the complex twisted homology of a nearby smooth curve splits into the homology groups of the component regions $C_0,C_l,C_f,C_{r-1},\ldots,C_1$. Their dimensions are $1,1,2,1,\ldots,1$. Each component class extends to a compact contour in nearby fibers. Its oscillatory integral has a convergent expansion
\begin{equation}\label{eq:gamma-pure-sector}
 I_j(t,Q)=B_j^{-\hbar}
             \bigl(I_{j,0}(Q)+O(t)\bigr),
\end{equation}
where $I_{j,0}$ is the regularized component integral with its induced flat section. The expansion is locally uniform in $Q$ and meromorphic in $\hbar$.
\end{lemma}
\begin{proof}
Cut between successive scales in the quotient $w$-sphere. The intervening annuli have monodromies
\[
 \mathsf q^4,\quad\mathsf q^2,\quad
 \mathsf q^{-2},\ldots,\mathsf q^{-2(r-1)}.
\]
Assume $\mathsf q^{2m}\ne1$ for $1\le m\le r$. This also makes the monodromies $-\mathsf q^2$ and $\mathsf q^{2r}$ at the global pole circles nontrivial. An annulus with monodromy $\xi\ne1$ has twisted chain complex $\C\xrightarrow{\xi-1}\C$, hence vanishing $H_0$ and $H_1$. Mayer--Vietoris gives the direct sum over the component regions. The relative homology sequence identifies the absolute group with the group allowing ends at the global pole circles. A component region has three punctures, except $C_f$, which has four; their twisted first homology dimensions are respectively one and two.

The projections are obtained by cutting a contour at each annulus and closing its ends. For an oriented segment $\sigma$ with flat section, let $\ell_a,\ell_b$ be positive loops based at its initial and final ends, with inherited sections and monodromies $\xi_a,\xi_b$. Its compact regularization is
\begin{equation}\label{eq:gamma-cut-cap}
 \operatorname{reg}(\sigma)
 =\sigma+\frac{\ell_a}{\xi_a-1}-\frac{\ell_b}{\xi_b-1}.
\end{equation}
The loop boundaries cancel those of $\sigma$. The closing chains from the two sides of an annulus cancel after transport, and each can be moved to a fixed circle in the adjacent rescaled chart. Thus every component class has a compact representative in that chart. These operations take place over $\C(\mathsf q)$; the integral lattice remains the one defined in Section~\ref{sec:cycles}.

On a neighborhood of the compact representative, the factorization in Proposition~\ref{prop:semistable} writes $\lambda=t^{-b_j}$ times a nonvanishing factor holomorphic in $t$. After division by $B_j$, this factor tends to the component function with the prescribed branch. The pulled-back differential $d\log\mu$ is holomorphic in $t$ and tends to its component differential; on $C_0$ one retains the double cover. Uniform Taylor expansion on each simplex of the compact chain proves \eqref{eq:gamma-pure-sector}. Formula \eqref{eq:gamma-cut-cap} agrees with an open integral wherever its end contributions converge, and gives its meromorphic continuation elsewhere.
\end{proof}

\subsection{Gamma normalization on the components}\label{subsec:local-gamma}
On the negative arc $e_-$, put $\mu=-s$, $0<s<1$, and $x=-w=s+s^{-1}$. Along its orientation $x$ decreases from $\infty$ to $2$, while $d\mu/\mu=-dx/\sqrt{x^2-4}$. The branch on $e_-$ contributes $\mathsf q^{-r}$; multiplication by $\mathsf q^{r-1}$ in $N_0$ leaves $\mathsf q^{-1}$. Hence, for $\Re\hbar>0$,
\begin{equation}\label{eq:structure-cycle-real-integral}
 \int_{N_0}e^{-W/z}\frac{d\mu}{\mu}
  =\mathsf q^{-1}K^{-\hbar}
       \int_2^\infty (x^2-4)^{\hbar-1/2}
              \prod_{i=1}^l(x+c_i)^{-\hbar}\,dx.
\end{equation}
On the component $C_k$, put $x=c_ky$ and $m=r-k$, and cap the node ends of the segment inherited from $N_0$. Removing the unit factor $B_k^{-\hbar}$, the differential in \eqref{eq:structure-cycle-real-integral} tends on its compact representative to $\mathsf q^{-1}y^{-m\hbar-1}(1+y)^{-\hbar}\,dy$. Its period is therefore
\begin{equation}\label{eq:point-gamma-period}
 I_m=\mathsf q^{-1}\int_0^\infty{}^{\rm reg}
       y^{-m\hbar-1}(1+y)^{-\hbar}\,dy
     =\mathsf q^{-1}\mathrm B(-m\hbar,(m+1)\hbar).
\end{equation}
On $C_l$, the substitution $x=c_ly$ and removal of $B_L^{-\hbar}$ instead give $\mathsf q^{-1}y^{2\hbar-1}(1+y)^{-\hbar}\,dy$, with period $\mathsf q^{-1}\mathrm B(2\hbar,-\hbar)$. On $C_0$, removal of $B_R^{-\hbar}$ leaves
\[
 \mathsf q^{-1}\int_2^\infty{}^{\rm reg}
      (x^2-4)^{\hbar-1/2}\,dx
 =\mathsf q^{-1}2^{2\hbar-1}\mathrm B(-\hbar,\hbar+\tfrac12)
 =\mathsf q^{-1}\mathrm B(-\hbar,2\hbar).
\]
The substitution $x=2/\sqrt v$ gives the first equality, and Gamma duplication gives the second. Thus both terminal components have period $I_1$. Each beta evaluation is first made with independently variable endpoint exponents in a common convergence domain, continued meromorphically, and then specialized to the displayed multiples of $\hbar$; the capped contours give this continuation.

At an isolated fixed point of arm distance $m$, the tangent weights divided by $z$ are $-m\hbar,(m+1)\hbar,-\hbar$. The surface-supported structure sheaf contributes
\begin{align}\label{eq:point-gamma-normalization}
 &\frac{(1-\mathsf q^{-2})\Gamma(1-m\hbar)
       \Gamma(1+(m+1)\hbar)\Gamma(1-\hbar)}
      {m(m+1)\hbar^3}\notag\\
 &\qquad=(1-\mathsf q^{-2})\Gamma(-\hbar)
       \Gamma(-m\hbar)\Gamma((m+1)\hbar)
       =-\frac{\pi iz}{\nu}I_m.
\end{align}
Here the denominator is $e_T(T_pY)/z^3$, and the first factor is the normal-line Koszul factor. The last equality follows from
\[
 (1-\mathsf q^{-2})\Gamma(\hbar)\Gamma(-\hbar)
       =-\frac{2\pi i}{\hbar\mathsf q}.
\]
A fiber of weight $k\nu$ multiplies both the Gamma charge and its cycle by $\mathsf q^{-k}$. The components $C_l$ and $C_0$ use $m=1$.

On the trivalent component $C_f$, take $0<Q<1$ and
\[
 \Lambda_f(s)=\frac{(s-1)(s-Q)}s.
\]
Let $\gamma$ run from $0$ to $\infty$ along $s=-y$, $0<y<\infty$, with $\arg\Lambda_f=\pi$. Let $\beta$ be a counterclockwise circle $|s|=R$, where $Q<R<1$, carrying the branch
\[
 \Lambda_f^{-\hbar}
       =\mathsf q^{-1}(1-s)^{-\hbar}(1-Q/s)^{-\hbar}.
\]
The monodromies around $0$ and $Q$ cancel, so $\beta$ is a closed twisted cycle; $\gamma$ is a relative cycle joining the two poles of $\Lambda_f$. Their capped classes will also be denoted by $\gamma$ and $\beta$.

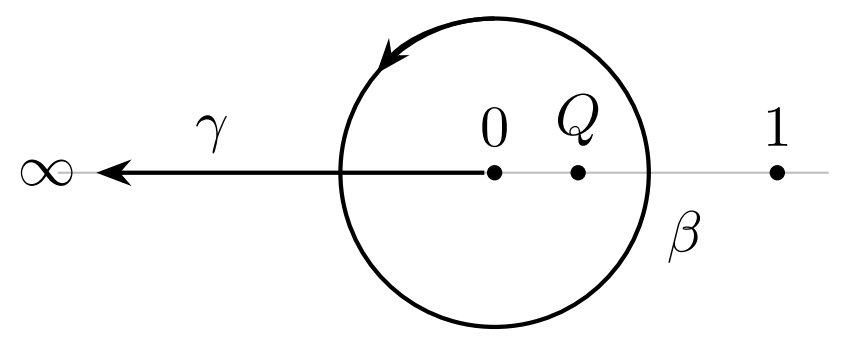
\begin{figure}[ht]
\centering
\begin{tikzpicture}[scale=.85,>=Stealth]
 \draw[gray!50] (-3.4,0)--(2.6,0);
 \draw[thick,->] (-.08,0)--(-3.1,0);
 \node[above] at (-2.2,0) {$\gamma$};
 \node[left] at (-3.1,0) {$\infty$};
 \draw[thick] (0,0) circle[radius=1.2];
 \draw[thick,->] (0,1.2) arc[start angle=90,end angle=140,radius=1.2];
 \node[right] at (1.18,-.5) {$\beta$};
 \foreach \x/\lab in {0/{0},.65/{Q},2.2/{1}} {
  \fill (\x,0) circle (1.7pt); \node[above] at (\x,.05) {$\lab$};
 }
\end{tikzpicture}
\caption{Cycles on the trivalent component $C_f$, with $0<Q<1$. The path $\gamma$ joins the poles $0$ and $\infty$ of $\Lambda_f$; $\beta$ surrounds its pole $0$ and zero $Q$.}
\label{fig:trivalent-cycles}
\end{figure}

Euler's integral and coefficient extraction give
\begin{align}\label{eq:fork-gamma-periods}
 I_\gamma(Q)&:=\int_\gamma\Lambda_f^{-\hbar}\frac{ds}s
 =\mathsf q^{-1}\mathrm B(\hbar,\hbar)
                   \,{}_2F_1(\hbar,\hbar;2\hbar;1-Q),\\
 I_\beta(Q)&:=\int_\beta\Lambda_f^{-\hbar}\frac{ds}s
 =2\pi i\,\mathsf q^{-1}\,{}_2F_1(\hbar,\hbar;1;Q).
 \notag
\end{align}
Indeed, the substitution $y=Qu/(1-u)$ in
\[
 I_\gamma=\mathsf q^{-1}\int_0^\infty
          \frac{y^{\hbar-1}\,dy}{(y+1)^\hbar(y+Q)^\hbar}
\]
gives the first formula for $\Re\hbar>0$. The second is $2\pi i\mathsf q^{-1}$ times the constant term of $(1-s)^{-\hbar}(1-Q/s)^{-\hbar}$. The formulas extend the periods meromorphically in $\hbar$.

Let $i_f$ be the surface zero section in the local threefold $\operatorname{Tot}(\OO_E(-2)_{2\nu}\oplus\OO_{E,-2\nu})$, and let $L$ be the pullback of $\OO_E(1)$ to the local surface, with trivial character along $E$. For a K-class $F$ on this surface, define $C_0(\mathsf q)$ and $C_1(\mathsf q)$ by
\[
 \widetilde{\operatorname{ch}}_z(F)|_E
      =C_0(\mathsf q)-\frac{2\pi i}{z}C_1(\mathsf q)h,
 \qquad h=c_1(L).
\]
Thus $F=\mathsf u^kL^d$ has $C_0=\mathsf q^{-k}$ and $C_1=d\mathsf q^{-k}$.

\begin{lemma}\label{lem:local-gamma}
At unit coordinate zero, the $T$-equivariant Gamma central charge of the trivalent theory satisfies
\begin{equation}\label{eq:fork-complete-gamma}
 Z_{(i_f)_*F}=-\frac{\pi iz}{\nu}
                  (C_0 I_\gamma-C_1 I_\beta).
\end{equation}
In particular, its cycle assignment is
\begin{equation}\label{eq:fork-k-cycles}
 (i_f)_*\mathsf u^kL^d
       \longmapsto\mathsf q^{-k}(\gamma-d\beta).
\end{equation}
\end{lemma}
\begin{proof}
The tangent roots along $E$ are $2h,2\nu-2h,-2\nu$. Since $h^2=0$, their Gamma product is
\[
 \Gamma(1+\hbar)\Gamma(1-\hbar)
 \left(1+\frac{2h}{z}
              \bigl(\psi(1)-\psi(1+\hbar)\bigr)\right).
\]
Multiplying by the Koszul factor $1-\mathsf q^{-2}$ and the modified Chern character, then pairing with $J^{\rm norm}(-z)=1-h\log Q/z+O(Q\log Q)$ using \eqref{eq:curve-metric}, gives
\begin{align}\label{eq:fork-gamma-asymptotic}
 Z_{(i_f)_*F}={}&-(1-\mathsf q^{-2})\Gamma(\hbar)\Gamma(-\hbar)\nonumber\\
 &\quad\cdot\bigl[C_0(\log Q+2\gamma_E+2\psi(\hbar))
                         +2\pi i C_1\bigr]+O(Q\log Q).
\end{align}
Here $\psi=\Gamma'/\Gamma$ and $\gamma_E=-\psi(1)$; we used $\psi(1+\hbar)=\psi(\hbar)+1/\hbar$.

The zero-balanced hypergeometric expansion in \eqref{eq:fork-gamma-periods} gives
\[
 I_\gamma=\mathsf q^{-1}
       [-\log Q-2\gamma_E-2\psi(\hbar)+O(Q\log Q)],
 \qquad
 I_\beta=2\pi i\mathsf q^{-1}+O(Q).
\]
Together with $(1-\mathsf q^{-2})\Gamma(\hbar)\Gamma(-\hbar)\mathsf q=-\pi iz/\nu$, these are precisely the logarithmic coefficient and constant term of \eqref{eq:fork-gamma-asymptotic}.

Both periods satisfy
\begin{equation}\label{eq:fork-gamma-qde}
 D^2I=\frac{2\hbar Q}{1-Q}DI
            +\frac{\hbar^2Q}{1-Q}I,
 \qquad D=Q\partial_Q.
\end{equation}
The scalar quantum equation obtained from \eqref{eq:curve-product} at argument $-z$ is the same equation, so it is also satisfied by the Gamma central charge. Its recurrence determines a solution from the coefficient of $\log Q$ and the constant term at $Q=0$: for the analytic solution $\sum_{n\ge0}a_nQ^n$, it is $n^2a_n=(n-1+\hbar)^2a_{n-1}$, and the logarithmic solution is then determined up to that analytic solution. This proves \eqref{eq:fork-complete-gamma} near $Q=0$. Analytic continuation in $Q$ and meromorphic continuation in $\hbar$ prove the asserted identity. A unit coordinate $\tau_0$ multiplies both sides by $e^{-\tau_0/z}$.
\end{proof}

The local periods also satisfy
\begin{equation}\label{eq:fork-gamma-wronskian}
 I_\gamma D I_\beta-I_\beta D I_\gamma
       =2\pi i\mathsf q^{-2}(1-Q)^{-2\hbar}.
\end{equation}
This follows from \eqref{eq:fork-gamma-qde} and the expansions at $Q=0$. Thus a trivalent period and its $Q$-derivative recover its two coefficients in the basis $I_\gamma,I_\beta$.

\subsection{The line-bundle contours on the components}\label{subsec:gamma-projections}
The local comparisons specify the component cycles for each fixed-locus restriction of a line bundle. We now verify that the global integral contours of Section~\ref{sec:cycles} have precisely these projections.
Write $N_0,N_a,N_L,N_R$ for the contours in \eqref{eq:line-bundle-cycles}. In a point column of the following table, a scalar denotes that multiple of the corresponding component of $N_0$. The trivalent column uses the cycles just defined.
\begin{proposition}\label{prop:gamma-contour-projections}
The projections of the fundamental contours under Lemma~\ref{lem:gamma-component-splitting} are
\begin{equation}\label{eq:gamma-contour-projections}
\begin{array}{c|c|c|c|c}
 &C_k,\ 1\le k<r&C_f&C_l&C_0\\ \hline
 N_0&1&\gamma&1&1\\
 N_a,\ a<r
 &\begin{cases}\mathsf q^{-2a},&a<k,\\
                 \mathsf q^{-2r},&a\ge k\end{cases}
 &\mathsf q^{-2a}\gamma&\mathsf q^{-2a}&\mathsf q^{-2a}\\
 N_r&\mathsf q^{-2r}&\mathsf q^{-2r}(\gamma+\beta)
             &\mathsf q^{-2r}&\mathsf q^{-2r}\\
 N_L&\mathsf q^{-r}&\mathsf q^{-r}\gamma
             &\mathsf q^{-r-2}&\mathsf q^{-r}\\
 N_R&\mathsf q^{-r}&\mathsf q^{-r}\gamma
             &\mathsf q^{-r}&\mathsf q^{-r-2}
\end{array}
\end{equation}
\end{proposition}
\begin{proof}
On $C_k$ use $\mu=u/\kappa_k$ and put $m=r-k$. After removing $B_k^{-\hbar}$, the differential on the positive segment before $u=1$ tends to
\[
 \mathsf q^{1-k}u^{(m+1)\hbar-1}(1-u)^{-\hbar}\,du.
\]
Let
\[
 G=\mathrm B((m+1)\hbar,-m\hbar),\quad
 H=\mathrm B((m+1)\hbar,1-\hbar),\quad
 J=\mathrm B(-m\hbar,1-\hbar).
\]
An arc ending at the $k$th zero contributes $\mathsf q^{1-k}H$. An arc ending later, including $e_+$, contributes $\mathsf q^{1-k}(H+\mathsf q^{-1}J)$, since passage below the zero changes the branch by $\mathsf q^{-1}$. The negative arc $e_-$ contributes $\mathsf q^{-r}G$. Reflection gives
\[
 \frac HG=-[m]_{\mathsf q},\qquad
 \frac JG=[m+1]_{\mathsf q},\qquad
 -[m]_{\mathsf q}+\mathsf q^{-1}[m+1]_{\mathsf q}
     =\mathsf q^{-m-1},\qquad
 [m]_{\mathsf q}=\frac{\mathsf q^m-\mathsf q^{-m}}
                         {\mathsf q-\mathsf q^{-1}}.
\]
The component has one-dimensional homology and its $N_0$ period is the nonzero meromorphic function $\mathsf q^{-1}G$. Over $\C(\mathsf q)$, the zero-end regularization is
$\operatorname{reg}(\gamma_j)=\gamma_j+\mathsf q^2f^{-1}\ell_j$,
so $\delta_j=f\,\operatorname{reg}(\gamma_j)$ in the pole-relative group. Dividing the regularized arc periods by the $N_0$ period therefore determines the component classes:
\begin{equation}\label{eq:point-arc-projections}
\begin{aligned}
 e_\pm&\longmapsto\mathsf q^{1-r},\\
 \delta_j&\longmapsto
 \begin{cases}
 0,&j<k-1,\\
 \mathsf q^{3-r}-\mathsf q^{r+3-2k},&j=k-1,\\
 (\mathsf q^2-1)\mathsf q^{1-r},&j\ge k.
 \end{cases}
\end{aligned}
\end{equation}
An earlier arc ends outside this component region and has zero projection. The beta identities are identities of compact regularized periods, obtained first with independent end exponents and then continued meromorphically.

On $C_l$, the positive differential is $\mathsf q^{-r-1}u^{-\hbar-1}(1-u)^{-\hbar}du$. Put $G=\mathrm B(-\hbar,2\hbar)$. The last zero arc contributes $\mathsf q^{-r-1}\mathrm B(-\hbar,1-\hbar)$; passage through that zero adds $\mathsf q^{-r-2}\mathrm B(2\hbar,1-\hbar)$. The ratios of these beta values to $G$ are $\mathsf q+\mathsf q^{-1}$ and $-1$. Hence
\[
 e_\pm\longmapsto\mathsf q^{1-r},\qquad
 \delta_j\longmapsto0\ (j<r+1),\qquad
 \delta_{r+1}\longmapsto\mathsf q^{3-r}-\mathsf q^{-r-1}.
\]
On $C_0$, all zero arcs disappear. Writing $\mu=-x$ on $e_-$ gives, after removal of $B_R^{-\hbar}$,
\[
 \mathsf q^{-r}\int_0^1{}^{\rm reg}
 x^{-2\hbar-1}(1-x^2)^{2\hbar}dx
 =\frac{\mathsf q^{-r}}2\mathrm B(-\hbar,1+2\hbar)
 =\mathsf q^{-r}G.
\]
The positive arc has crossed all $r+2$ zeros, so its coefficient is $\mathsf q^{-r-2}$. Therefore on this component $e_-\mapsto\mathsf q^{1-r}$ and $e_+\mapsto\mathsf q^{-r-1}$.

On $C_f$, use $\mu=u/\kappa_r$. The positive differential before its first zero is
\[
 \mathsf q^{1-r}u^{\hbar-1}
            (1-u)^{-\hbar}(1-Qu)^{-\hbar}\,du.
\]
The arc ending at $u=1$ has period $\mathsf q^{1-r}\mathrm B(\hbar,1-\hbar)\,{}_2F_1(\hbar,\hbar;1;Q)$. Comparing with $I_\beta$ gives $\delta_{r-1}\mapsto\mathsf q^{3-r}\beta$. The through-going positive arc can be deformed below both zeros to the negative ray. For $0<\Re\hbar<1$ the end contributions vanish and this rotation supplies $\mathsf q^{-1}$, giving the same period as $e_-$. Thus
\[
 e_\pm\longmapsto\mathsf q^{1-r}\gamma,\qquad
 \delta_{r+1}\longmapsto(\mathsf q^2-1)\mathsf q^{1-r}\gamma,
 \qquad \delta_j\longmapsto0\quad(j<r-1).
\]
These are homology identities: the deformation takes place in the punctured component, with the continued flat sections.

To compute the remaining contour, let $\beta_t$ be the counterclockwise circle $|w|=Rc_r$, with $Q<R<1$, carrying the section continued from $e_-$ at $w=-Rc_r$. It surrounds $c_{r+1},c_l,2,-2$, whose monodromies have product one. Cutting this circle along the arcs gives
\begin{equation}\label{eq:global-beta-relation}
 [\beta_t]=\delta_r+\mathsf q^2\delta_{r+1}
              -(1+\mathsf q^2)(\mathsf q^2e_+-e_-).
\end{equation}
In the chart $s=w/c_r$, the lift near $\mu=0$ has $d\log\mu\to-ds/s$. The section inherited from $e_-$ and the scalar $B_f^{-\hbar}$ therefore identify $-\mathsf q^{r-1}\beta_t$ with $\beta$. Substitution in \eqref{eq:global-beta-relation} gives
\[
 \delta_r\longmapsto
 \mathsf q^{1-r}\bigl((\mathsf q^2-1)\gamma-\beta\bigr).
\]
Finally insert these arc projections into \eqref{eq:line-bundle-cycles}. The finite sums telescope; for the long-arm rows one may use $\mathsf q^2N_{a+1}-N_a=\mathsf q^{-r-1}\delta_a$ for $a<r$. This gives every entry of \eqref{eq:gamma-contour-projections}, including $r=2$.
\end{proof}

\subsection{Determination by the boundary data}\label{subsec:gamma-uniqueness}
The scalar equation \eqref{eq:scalar-period-qde} makes $v=(\partial_a I)_a$ a flat covector of the quantum connection. Along the family \eqref{eq:arc} its equations have the form
\begin{equation}\label{eq:gamma-normal-system}
 t\partial_t v=A(t,Q)v,\qquad \partial_Qv=B(t,Q)v,
\end{equation}
with $A,B$ holomorphic at $t=0$ for $Q\ne0,1$. Indeed the root expression for the quantum product is holomorphic in the vanishing Novikov variables at this face. Each transverse divisor has degree zero on $E$, so the residue acts by a scalar on each point block and on the whole rank-two curve block. In the fixed-component frame,
\[
 A(0,Q)=A_0
 =\hbar\,\operatorname{diag}(b_1,\ldots,b_{r-1},b_f,b_f,b_L,b_R),
\]
with \eqref{eq:gamma-normal-exponents}. These eigenvalues are also obtained by restricting $\tau$ to the fixed components:
\begin{equation}\label{eq:gamma-unit-restrictions}
 \tau|_{p_j}=2\nu\log B_j,\qquad
 \tau|_E=2\nu\log B_f+h\log Q.
\end{equation}
\rev{Here $h=c_1(\OO_E(1))$, and the equalities follow from \eqref{eq:mirror-map} and $c_1^T(\mathcal L_i)=-D_i$.} The factors $a_j(Q)^{-\hbar}$ in $B_j^{-\hbar}$ retain the unit shifts when $Q$ varies.

\begin{lemma}\label{lem:gamma-boundary-uniqueness}
Suppose $(b_i-b_j)\hbar\notin\mathbb Z$ for distinct component exponents. A solution of \eqref{eq:scalar-period-qde} is determined by the leading scalar coefficient in each point sector and the complete leading scalar function of $Q$ in the curve sector.
\end{lemma}
\begin{proof}
There is a unique convergent matrix $H(t,Q)=I+O(t)$ satisfying
\[
 t\partial_tH=AH-HA_0.
\]
Writing $A=A_0+\sum_{n\ge1}A_nt^n$ gives the recursion
\[
 (nI-\operatorname{ad}A_0)H_n
             =\sum_{k=1}^n A_kH_{n-k}.
\]
The operators on the left are invertible, with inverses of norm $O(n^{-1})$, so a power-series majorant proves convergence locally uniformly in $Q$. After this gauge change, compatibility of \eqref{eq:gamma-normal-system} gives $t\partial_t\widetilde B=[A_0,\widetilde B]$. The same nonresonance makes every positive-degree coefficient vanish, hence $\widetilde B=B(0,Q)$. Every flat covector therefore has the form
\begin{equation}\label{eq:gamma-frobenius-solution}
 v=H(t,Q)t^{A_0}v_0(Q),\qquad
 \partial_Qv_0=B(0,Q)v_0.
\end{equation}

On a point block the unit is one, so $I=-z\partial_0I$ recovers its one scalar coordinate. On the curve block, remove the unit factor $a_f(Q)^{-\hbar}$. The scalar function and its $Q$-derivative give the two coordinates of the local quantum solution: the unit and $h$ span this block, and \eqref{eq:fork-gamma-wronskian} is nonzero. Thus the specified scalar data determine all of $v_0(Q)$. Formula \eqref{eq:gamma-frobenius-solution} determines the full gradient at an interior point, and the flat quantum connection continues it uniquely in all the parameters. The unit equation recovers $I$ itself.
\end{proof}

\begin{proof}[Proof of Theorem~\ref{thm:gamma-period}]
First take generic $\hbar$ satisfying the monodromy and nonresonance conditions above. By Proposition~\ref{prop:scalar-period-qde}, the two sides satisfy the same scalar quantum equations. Lemma~\ref{lem:gamma-component-splitting} expresses each contour as a sum of compact component contours, whose periods have the characteristic-power expansions \eqref{eq:gamma-pure-sector}. Since distinct exponents differ by nonintegers, these are precisely the separate sectors of \eqref{eq:gamma-frobenius-solution}.

On the Gromov--Witten side, Proposition~\ref{prop:split} specializes the quantum theory to the point and curve summands. Its normalized $J$-function specializes in each sector to the corresponding local $J$-function, with the divisor exponentials \eqref{eq:gamma-unit-restrictions}; this follows coefficientwise from localization of the one-point genus-zero series, including its unstable completion. Terms of positive degree in a transverse Novikov variable give positive powers of $t$. Thus its boundary coefficients are the local Gamma charges of the restrictions of the K-class. Equivalently, these coefficients specify the formal solution \eqref{eq:gamma-frobenius-solution}, whose convergence follows from the same recursion.
In the fixed-component frame this specialization reads
\[
 \begin{aligned}
 J_Y^{\rm norm}(\tau(t,Q),-z)|_{p_j}
   &=B_j^{-\hbar}(1+O(t)),\\
 J_Y^{\rm norm}(\tau(t,Q),-z)|_E
   &=B_f^{-\hbar}
       \bigl(J_{\rm loc}^{\rm norm}(Q,-z)+O(t)\bigr),
 \end{aligned}
\]
where $J_{\rm loc}^{\rm norm}$ is the normalized $(-2,0)$ curve $J$-function at unit coordinate zero. Thus the curve coefficient retains its full $Q$-dependence.

For the fundamental line bundles, Lemma~\ref{lem:bundle-restrictions} and Proposition~\ref{prop:gamma-contour-projections} identify these restrictions with their component contours. At a point, the character $\mathsf u^k$ becomes $\mathsf q^{-k}$ and \eqref{eq:point-gamma-normalization} supplies the factor $-\pi iz/\nu$. On $E$, the restriction $\mathsf u^k\OO_E(d)$ corresponds to $\mathsf q^{-k}(\gamma-d\beta)$ by Lemma~\ref{lem:local-gamma}, with the same factor. Every point coefficient and the full rank-two curve solution therefore agree. Lemma~\ref{lem:gamma-boundary-uniqueness} proves the identity near the partial boundary and then throughout the chosen large-radius domain.

The fundamental supported line bundles form an integral basis by Corollary~\ref{cor:line-bundle-cycles}. Additivity and character twists prove the result for every $F\in K_Z^T(Y)$. Meromorphic continuation in $\hbar$ removes the genericity restrictions, and transporting the contours with their flat sections gives the stated continuation in the mirror parameters.
\end{proof}

\subsection{Gamma descendants}\label{sec:gamma-descendants}
Section~\ref{sec:correlators} has already established the all-genus descendant dictionary. It remains to realize its leaf substitution by integration over the cycles constructed in Section~\ref{sec:cycles}.
We use the integration-by-parts argument for descendant leaves in \cite[Lemma~7.3]{FLYZ}, with the normalization of our sign-sector recursion. For a supported class $F$, write $\gamma_F:=\operatorname{mir}_Z(F)\in\mathcal H_{\mathscr R}^+$ for its cycle on the smooth global mirror curve, as defined in Section~\ref{sec:cycles}. In the conventions of Theorem~\ref{thm:gamma-period},
\[
 Z_F=-\frac{\pi iz}{\nu}
       \int_{\gamma_F}^{\rm reg}e^{-W/z}\,d\log\mu.
\]
Set $\rho_\alpha=\sqrt2s_\alpha$ and
\[
 \mathcal B_\alpha(p)=
     \left.\frac{B^-(p,q)}{d\rho_\alpha(q)}\right|_{q=p_\alpha},
 \qquad \omega_a=(\partial_aW)\sqrt\nu\,d\log\mu.
\]
Our square roots give $\thh_\alpha^0=i\mathcal B_\alpha$ and
\begin{equation}\label{eq:gamma-leaf-metric}
 \left.\frac{\omega_a}{d\rho_\alpha}\right|_{p_\alpha}
        =2\nu^{3/2}\eta(T_a,\eh_\alpha),
\end{equation}
where $T_a$ is the flat vector associated to $\partial_a$. On the rational cover one has
\begin{equation}\label{eq:gamma-cauchy}
 \sum_\alpha
 \left.\frac{\omega_a}{d\rho_\alpha}\right|_{p_\alpha}
       \mathcal B_\alpha
       =-d\left(\frac{\omega_a}{d\xh}\right).
\end{equation}
The ratio on the right is an anti-invariant rational function whose poles are the simple critical points. Its logarithmic puncture poles cancel. Both sides have the same double principal parts at the two lifts of every critical orbit; their difference is therefore a holomorphic one-form on $\mathbb P^1$ and vanishes.

The $T$-equivariant $J$-function and $S$-operator satisfy
\[
 \partial_aZ_F=-\frac1z
      \eta(T_a,S_\tau(-z)\kappa_z(F)).
\]
Differentiating the oscillatory integral along a flat family of cycles gives
\[
 \partial_aZ_F=\frac{\pi i}{\nu^{3/2}}
               \int_{\gamma_F}^{\rm reg}e^{-W/z}\omega_a.
\]
Integrating \eqref{eq:gamma-cauchy} by parts and using \eqref{eq:gamma-leaf-metric} now yields
\[
 2\pi\int_{\gamma_F}^{\rm reg}e^{-W/z}\thh_\alpha^0
 =-z^{-2}\eta(\eh_\alpha,S_\tau(-z)\kappa_z(F)).
\]
To derive the undressed recursion, put $\mathcal D(\xi)=-d(\xi/d\xh)$ and write $R(z)=I+\sum_{m\ge1}R_mz^m$. The definition \eqref{eq:theta}, symmetry of the kernel, and \eqref{eq:r-jets} give
\[
 \mathcal D(\theta_\alpha^k)
 =\theta_\alpha^{k+1}
       -\sum_\beta(R_{k+1})_{\beta\alpha}\theta_\beta^0.
\]
Indeed, the highest polar part at $p_\alpha$ gives the first term. At $p_\beta$ the regular constant of $\theta_\alpha^k$ is
$i(2k-1)!!\,b_{2k}^{\beta\alpha}d\rho_\beta$.
Dividing by $d\xh=\rho_\beta d\rho_\beta$ and differentiating produces the double pole recorded by the second term, since
$(R_{k+1})_{\beta\alpha}=-(2k-1)!!b_{2k}^{\beta\alpha}$.
There are no additional poles at the logarithmic punctures. The difference is therefore a holomorphic one-form on the rational cover and vanishes.

In row-vector notation \eqref{eq:leaf-dressing} is $\theta(z)=\thh(z)R(z)$. The preceding identity yields
\[
 \mathcal D\theta(z)=\frac{\theta(z)-\theta(0)R(z)}z.
\]
As $R$ is independent of the curve variable and invertible, this gives
$\mathcal D\thh(z)=(\thh(z)-\thh(0))/z$, or
\[
 \thh_\alpha^{k+1}=-d(\thh_\alpha^k/d\xh).
\]
Another integration by parts multiplies the oscillatory integral by $z^{-1}$. Thus
\begin{equation}\label{eq:gamma-leaf-integral}
 2\pi\int_{\gamma_F}^{\rm reg}e^{-W/z}\thh_\alpha^k
 =-z^{-k-2}\eta(\eh_\alpha,S_\tau(-z)\kappa_z(F)).
\end{equation}
These operations first hold in a decay sector. The resulting identities extend meromorphically. The weighted critical residues vanish: each $\thh_\alpha^k$ is a finite linear combination of the $\theta_\beta^j$, whose critical principal parts contain only even negative powers $s_\beta^{-2m}ds_\beta$. Since $e^{-W/z}=e^{\xh/z}$ is an even function of $s_\beta$, their product has no $s_\beta^{-1}ds_\beta$ term. The regularization therefore also allows deformation of the contours through the critical disks.

\begin{theorem}\label{thm:gamma-descendants}
For $F_1,\ldots,F_n\in K_Z^T(Y)$, with $n>0$ and $2g-2+n>0$, the $T$-equivariant descendant generating series satisfies
\begin{align}\label{eq:gamma-descendants}
 &\sum_{\beta,m}\frac{1}{m!}
 \left\langle
   \prod_{j=1}^n\frac{\kappa_{z_j}(F_j)}{z_j(z_j+\psi_j)},
          \tau,\ldots,\tau\right\rangle_{g,n+m,\beta}^{Y}
 \nonumber\\
 &\hspace{15mm}=(2\pi)^n
 \int_{\gamma_{F_1}\times\cdots\times\gamma_{F_n}}^{\rm reg}
 e^{-\sum_jW(p_j)/z_j}\,\omega_{g,n}(p_1,\ldots,p_n).
\end{align}
Each leg uses $\mathsf q=e^{2\pi i\nu/z_j}$ in its cycle coefficients. The identity is interpreted coefficientwise in the descendant expansions and then continued along the chosen smooth chamber. Equivalently, write the left-hand side with Novikov weight $q^\beta$, extra primary insertions $t_0\,1$, and leg classes $e^{-(\tau-t_0\,1)/z_j}\kappa_{z_j}(F_j)$, where $\tau-t_0\,1=\sum_i(\log q_i)D_i$ acts by cup product. This follows from the descendant divisor equation, since for $u(s)=\kappa_z(F)/(z(z+s))$ its correction is $D_i(u(s)-u(0))/s=-D_i u(s)/z$.
\end{theorem}
\begin{proof}
For the input $u(s)=\kappa_z(F)/(z(z+s))$, the ancestor transform has coefficients
\[
 U_k^\alpha=(-1)^kz^{-k-2}
            \eta(\eh_\alpha,S_\tau(-z)\kappa_z(F)).
\]
Equation~\eqref{eq:gamma-leaf-integral} is precisely the substitution $\thh_\alpha^k\mapsto(-1)^{k+1}U_k^\alpha$. Apply \eqref{eq:descendant-dictionary} separately to every leg. At fixed genus, number of markings, and descendant coefficient, only finitely many terms of the stable graph expansion occur.
\end{proof}


\begin{thebibliography}{99}
\bibitem{AK}
K.~Aomoto and M.~Kita, \emph{Theory of Hypergeometric Functions}, Springer Monographs in Mathematics, Springer, Tokyo, 2011.
\bibitem{Behrend}
K.~Behrend, \emph{Gromov--Witten invariants in algebraic geometry}, Invent.\ Math.\ \textbf{127} (1997), 601--617.
\bibitem{BF}
K.~Behrend and B.~Fantechi, \emph{The intrinsic normal cone}, Invent.\ Math.\ \textbf{128} (1997), 45--88.
\bibitem{BB}
G.~Borot and A.~Brini, \emph{Chern--Simons theory on spherical Seifert manifolds, topological strings and integrable systems}, Adv.\ Theor.\ Math.\ Phys.\ \textbf{22} (2018), no.~2, 305--394.
\bibitem{BKMP}
V.~Bouchard, A.~Klemm, M.~Mari\~no, and S.~Pasquetti, \emph{Remodeling the B-model}, Comm.\ Math.\ Phys.\ \textbf{287} (2009), 117--178.
\bibitem{BKR}
T.~Bridgeland, A.~King, and M.~Reid, \emph{The McKay correspondence as an equivalence of derived categories}, J.\ Amer.\ Math.\ Soc.\ \textbf{14} (2001), 535--554.
\bibitem{BriniE8}
A.~Brini, \emph{$E_8$ spectral curves}, Proc.\ Lond.\ Math.\ Soc.\ (3) \textbf{121} (2020), no.~4, 954--1032.
\bibitem{BM}
A.~Brini and J.~Ma, \emph{Integral and higher-genus quantum McKay correspondence in type D}, arXiv:2610.22359 (2026).
\bibitem{BMS}
A.~Brini, J.~Ma, and I.~A.~B.~Strachan, \emph{Dubrovin duality and mirror symmetry for ADE resolutions}, Proc.\ Roy.\ Soc.\ A \textbf{481} (2025), article 20250047.
\bibitem{BvG}
A.~Brini and K.~van Gemst, \emph{Mirror symmetry for extended affine Weyl groups}, J.\ \'{E}c.\ polytech.\ Math.\ \textbf{9} (2022), 907--957.
\bibitem{BG}
J.~Bryan and A.~Gholampour, \emph{Root systems and the quantum cohomology of ADE resolutions}, Algebra Number Theory \textbf{2} (2008), 369--390.
\bibitem{CCIT}
T.~Coates, A.~Corti, H.~Iritani, and H.-H.~Tseng, \emph{A mirror theorem for toric stacks}, Compos.\ Math.\ \textbf{151} (2015), 1878--1912.
\bibitem{CG}
T.~Coates and A.~Givental, \emph{Quantum Riemann--Roch, Lefschetz and Serre}, Ann.\ of Math.\ \textbf{165} (2007), 15--53.
\bibitem{DSZZ}
B.~Dubrovin, I.~A.~B.~Strachan, Y.~Zhang, and D.~Zuo, \emph{Extended affine Weyl groups of BCD type, Frobenius manifolds and their Landau--Ginzburg superpotentials}, Adv.\ Math.\ \textbf{351} (2019), 897--946.
\bibitem{DOSS}
P.~Dunin-Barkowski, N.~Orantin, S.~Shadrin, and L.~Spitz, \emph{Identification of the Givental formula with the spectral curve topological recursion procedure}, Comm.\ Math.\ Phys.\ \textbf{328} (2014), 669--700.
\bibitem{Eynard}
B.~Eynard, \emph{Invariants of spectral curves and intersection theory of moduli spaces of complex curves}, Commun.\ Number Theory Phys.\ \textbf{8} (2014), 541--588.
\bibitem{EO}
B.~Eynard and N.~Orantin, \emph{Invariants of algebraic curves and topological expansion}, Commun.\ Number Theory Phys.\ \textbf{1} (2007), 347--452.
\bibitem{EOBKMP}
B.~Eynard and N.~Orantin, \emph{Computation of open Gromov--Witten invariants for toric Calabi--Yau 3-folds by topological recursion, a proof of the BKMP conjecture}, Comm.\ Math.\ Phys.\ \textbf{337} (2015), 483--567.
\bibitem{Fang}
B.~Fang, \emph{Central charges of T-dual branes for toric varieties}, Trans.\ Amer.\ Math.\ Soc.\ \textbf{373} (2020), 3829--3851.
\bibitem{FLM}
B.~Fang, Z.~Lan, and J.~Ma, \emph{Involution-equivariant topological recursion and mirror symmetry for the affine binary dihedral Calabi--Yau threefold}, \rev{arXiv:2607.07355v2} (2026).
\bibitem{FLYZ}
B.~Fang, C.-C.~M.~Liu, S.~Yu, and Z.~Zong, \emph{Remodeling Conjecture with Descendants}, arXiv:2504.15696 (2025).
\bibitem{FLZP1}
B.~Fang, C.-C.~M.~Liu, and Z.~Zong, \emph{The Eynard--Orantin recursion and equivariant mirror symmetry for the projective line}, Geom.\ Topol.\ \textbf{21} (2017), 2049--2092.
\bibitem{FLZ}
B.~Fang, C.-C.~M.~Liu, and Z.~Zong, \emph{On the remodeling conjecture for toric Calabi--Yau 3-orbifolds}, J.\ Amer.\ Math.\ Soc.\ \textbf{33} (2020), 135--222.
\bibitem{GKL}
A.~Giacchetto, R.~Kramer, and D.~Lewa\'nski, \emph{A new spin on Hurwitz theory and ELSV via theta characteristics}, Selecta Math.\ (N.S.)\ \textbf{31} (2025), article 90.
\bibitem{GiventalHigher}
A.~B.~Givental, \emph{Semisimple Frobenius structures at higher genus}, Internat.\ Math.\ Res.\ Notices (2001), no.~23, 1265--1286.
\bibitem{GiventalQuant}
A.~B.~Givental, \emph{Gromov--Witten invariants and quantization of quadratic Hamiltonians}, Mosc.\ Math.\ J.\ \textbf{1} (2001), 551--568.
\bibitem{GSV}
G.~Gonzalez-Sprinberg and J.-L.~Verdier, \emph{Construction g\'eom\'etrique de la correspondance de McKay}, Ann.\ Sci.\ \'Ecole Norm.\ Sup.\ (4) \textbf{16} (1983), 409--449.
\bibitem{GP}
T.~Graber and R.~Pandharipande, \emph{Localization of virtual classes}, Invent.\ Math.\ \textbf{135} (1999), 487--518.
\bibitem{Hosono}
S.~Hosono, \emph{Central charges, symplectic forms, and hypergeometric series in local mirror symmetry}, in \emph{Mirror Symmetry V}, AMS/IP Stud.\ Adv.\ Math.\ \textbf{38}, Amer.\ Math.\ Soc., Providence, RI, 2006, 405--439.
\bibitem{Iritani}
H.~Iritani, \emph{An integral structure in quantum cohomology and mirror symmetry for toric orbifolds}, Adv.\ Math.\ \textbf{222} (2009), 1016--1079.
\bibitem{IN}
Y.~Ito and I.~Nakamura, \emph{Hilbert schemes and simple singularities}, in \emph{New Trends in Algebraic Geometry}, London Math.\ Soc.\ Lecture Note Ser.\ \textbf{264}, Cambridge University Press, Cambridge, 1999, 151--234.
\bibitem{KV}
M.~Kapranov and E.~Vasserot, \emph{Kleinian singularities, derived categories and Hall algebras}, Math.\ Ann.\ \textbf{316} (2000), 565--576.
\bibitem{KKP}
L.~Katzarkov, M.~Kontsevich, and T.~Pantev, \emph{Hodge theoretic aspects of mirror symmetry}, in \emph{From Hodge Theory to Integrability and TQFT: $tt^*$-Geometry}, Proc.\ Sympos.\ Pure Math.\ \textbf{78}, Amer.\ Math.\ Soc., Providence, RI, 2008, 87--174.
\bibitem{KM}
M.~Kontsevich and Y.~Manin, \emph{Relations between the correlators of the topological sigma-model coupled to gravity}, Comm.\ Math.\ Phys.\ \textbf{196} (1998), 385--398.
\bibitem{KrM}
O.~Kruglinskaya and A.~Marshakov, \emph{On Lie groups and Toda lattices}, J.\ Phys.\ A: Math.\ Theor.\ \textbf{48} (2015), article 125201.
\bibitem{MOPPZ}
A.~Marian, D.~Oprea, R.~Pandharipande, A.~Pixton, and D.~Zvonkine, \emph{The Chern character of the Verlinde bundle over $\overline{\mathcal M}_{g,n}$}, J.\ Reine Angew.\ Math.\ \textbf{732} (2017), 147--163.
\bibitem{MaThesis}
J.~Ma, \emph{Mirror symmetry and the quantum McKay correspondence}, Ph.D. thesis, University of Sheffield, 2026. \rev{\url{https://etheses.whiterose.ac.uk/id/eprint/39150/}}.
\bibitem{Marino}
M.~Mari\~no, \emph{Open string amplitudes and large order behavior in topological string theory}, J.\ High Energy Phys.\ 03 (2008), article 060.
\bibitem{Teleman}
C.~Teleman, \emph{The structure of 2D semi-simple field theories}, Invent.\ Math.\ \textbf{188} (2012), 525--588.
\bibitem{Thomason}
R.~W.~Thomason, \emph{Algebraic K-theory of group scheme actions}, in \emph{Algebraic Topology and Algebraic K-Theory}, Ann.\ of Math.\ Stud.\ \textbf{113}, Princeton University Press, Princeton, NJ, 1987, 539--563.
\bibitem{ThomasonLefschetz}
R.~W.~Thomason, \emph{Une formule de Lefschetz en K-th\'eorie \'equivariante alg\'ebrique}, Duke Math.\ J.\ \textbf{68} (1992), no.~3, 447--462.
\end{thebibliography}
\end{document}